\documentclass[12pt]{article}
\usepackage{graphicx}
\usepackage{latexsym,amssymb,amsmath,amscd,amsfonts,mathabx}
\topmargin - 2cm
\bibliography{alpha.bst}

\newtheorem{thm}{Theorem}[section]

\newtheorem{lem}[thm]{Lemma}

\newtheorem{definition}[thm]{Definition}
\newtheorem{rem}[thm]{Remark}
\newtheorem{ex}[thm]{Example}
\newcommand{\eproof}{\rule{0,2cm}{0,2cm}}

\newcommand{\vf}{\varphi}

\newcommand{\sing}{\mbox{sing}}

\newcommand{\re}{\mathbb{R}}

\begin{document}




\title{{\bf Representations of solutions of time-fractional multi-order systems of 
differential-operator equations }}

\author{Sabir Umarov$^{1}$, }
\date{}
\maketitle

\begin{center}
{\it $^1$University of New Haven, Department of Mathematics, \\300 Boston Post Road, West Haven, CT 06516, USA}
\end{center}
\vspace{10pt}

\begin{abstract} 
This paper is devoted to the general theory of systems of linear time-fractional differential-operator equations. The representation formulas for solutions of systems of ordinary differential equations with single (commensurate) fractional order is known through the matrix valued Mittag-Leffler function. Multi-order (incommensurate) systems with rational  components can be reduced to single order systems, and hence, representation formulas are also known. However, for arbitrary  fractional multi-order (not necessarily with rational components) systems of differential equations this question remains open even in the case of ordinary differential equations. In this paper we obtain representation formulas for solutions of arbitrary  fractional multi-order systems of differential-operator equations along with proving the existence and uniqueness theorems in appropriate topological-vector spaces. Moreover, we introduce vector-indexed Mittag-Leffler functions and prove some of their properties. 
\end{abstract}




\label{}
\section{Introduction}

Let $X$ be a reflexive Banach space and $A: {\mathcal{D}}
\rightarrow X$  a closed linear operator with a domain
${\mathcal{D}} \subset X.$  
Consider the systems of $m \ge 1$ time-fractional differential operator equations, which we will write in the form
\begin{equation}
 \label{01}
D^{\beta} \mathcal{U}(t) = \mathcal{F} (A) \mathcal{U}(t) + \mathcal{H}(t), \quad t>t_0,
\end{equation}
with the initial condition
\begin{equation}
\label{02}
B \, \mathcal{U}(t_0) = \varPhi.
\end{equation}
In equation \eqref{01} $t > t_0;$ $\mathcal{B}=(\beta_1,\dots,\beta_m), \  0 < \beta_j \le 1,$ fractional orders of the system;  
$
\mathcal{U}(t): [t_0,\infty) \to X \times \dots \times X,
$
is an abstract vector-valued function with components $u_j(t), \ j=1,\dots, m,$ to be found, and 
\begin{equation}
\label{D1}
D^{\mathcal{B}} \mathcal{U} = (D^{\beta_1}u_1(t), \dots, D^{\beta_m}u(t));
\end{equation}
Here $D^{\beta_j}, \ j=1,\dots, m,$ is the fractional order derivative of order $ 0<\beta_j \le 1$ in the sense of Riemann-Liouville or Caputo. The matrix-valued operator $\mathcal{F}(A)$ on the right hand side of equation \eqref{01} has the form
\begin{align}
\mathcal{F} (A)&= 
\begin{bmatrix}
f_{11} (A) & \dots & f_{1m} (A) \\
\dots && \\
f_{m1} (A) & \dots & f_{mm} (A)
\end{bmatrix};
\label{F0}
\end{align}
$H(t): [t_0,\infty] \to X \times \dots \times X $ is a given vector-valued function satisfying some conditions clarified later. 
In initial condition \eqref{02} the operator $B$ depends on whether $D^{\mathcal{B}}$ is in the sense of Riemann-Liouville or Caputo  and $\varPhi$ is given element of some topological-vector space specified later. The operator $\mathcal{F}(A)$ has a matrix-symbol $\mathcal{F}(z) \equiv \{f_{k i}(z)\}, k,j = 1, \dots, m, $  entries of which may have singularities in the spectrum of the operator $A.$ The exact definitions of operators $D^{\beta_j}, \ j=1,\dots,m,$ and $f_{k j} (A), \ k, j = 1,\dots,m,$ are given in Section \ref{sec_prel}. 

It is well known (see e.g. \cite{DalKrein,VasPis}) that if $\mathcal{B}=(1,\dots,1)$ then the solution is represented in the form 
\[
\mathcal{U}(t)=S(t,A)\varPhi + \int_{t_0}^t S(t-\tau,A) H(\tau) d\tau,
\]
where 
\[
S(t,A)=\exp(-t \mathcal{F}(A))
\] 
is {\it the solution operator,} which has an exponential form. It is also known \cite{SKM93,Ba2001,Kos93,Koc89,GLZ99,EidKoch2004,KilbasST,Umarov_book2015} that in the case of various fractional order differential equations (not systems) the solution can be represented through the Mittag-Leffler (ML) function 
\begin{equation}
\label{MLcl}
E_{\beta, \nu}(z) = \sum_{n=0}^{\infty} \frac{z^n}{\Gamma(\beta n +\nu)}, \quad z \in \mathbb{C},
\end{equation} 
which generalizes the exponential function. Namely, if $\beta= \nu=1,$ then we have $E_{\beta, \nu}(z)=\exp(z).$ In the case of systems, when $\mathcal{B}$ has equal components, that is $\beta_j=\beta, j=1,\dots,m,$ a representation formula for the solution is obtained via the matrix-valued ML function \cite{Veber,Garrapa,Varsha,UACh}.

In the case when the components of the vector-order $\mathcal{B}$ in equation \eqref{01}  are rational, that is $\beta_j = p_j/q_j,$ where $p_j$ and $q_j$ are coprime numbers,  the corresponding system can be reduced to a system with a scalar order $\alpha \in (0,1)$ \cite{DengLiGuo,good,Umarov2024}. However, the number of equations in the reduced system may increase significantly. Let $M$ be the lowest common multiple of numbers $q_1,\dots, q_m,$ and $\alpha = 1/M.$ Then, the number of equations in the reduced system becomes $N=M(\beta_1+\dots + \beta_m).$ For example, if the orders in the original system of 4 equations are $\beta_1=\frac{1}{2}, \beta_2= \frac{2}{3}, \beta_3= \frac{1}{5},$ and $\beta_4= \frac{6}{7},$ then $M=210$ and $N= 210 \cdot (1/2+2/3+1/5+6/7)=467.$   Thus, the reduced system will contain 467 equations of order $\alpha=\frac{1}{210},$ though originally we had only four equations in the system. Even numerical solutions of such reduced systems consume a significant computer and timing resources, and thus, the method of reducing to a scalar order systems should be considered ineffective. Therefore, developing the direct general techniques for solution and qualitative analysis of systems of fractional order differential equations with any positive real orders is important.

For a multi-order $\mathcal{B}$ with arbitrary components (not necessarily rational) representation formulas are found in \cite{Umarov2024} in some particular cases.  
In the current paper we obtain representation formulas for arbitrary multi-order $\mathcal{B}$ and matrix-valued operator $\mathcal{F}(A).$  We also introduce more general ML functions, called an {\it vector-indiced matrix-valued Mittag-Leffler function.} We show that the solution of system \eqref{01}-\eqref{02} is represented through an operator-dependent matrix-valued ML function. The results obtained in this paper are new even for time-fractional systems of linear ordinary differential equations.

We note that systems of fractional order ordinary and partial differential equations have rich applications. For example, they are used in modeling of processes in bio-systems \cite{DasGupta,Rihan,GuoFang}, ecology \cite{Khan,Rana}, epidemiology \cite{Zeb,Islam}, quantum systems \cite{TP,Ercan,Bernstein}, etc.

This paper is organized as follows. In Section \ref{sec_prel} we provide some preliminary facts on the Mittag-Leffler functions, including matrix valued versions. To our best knowledge Lemmas \ref{lemma_1} and \ref{lemma_2} are new. Here we also introduce the vector-indiced matrix-valued ML functions and study some of their properties used in this paper. This section familiarizes with the basic topological-vector spaces on which the corresponding matrix-valued operators with singular symbols act. In Sections \ref{sec_main} and \ref{sec_RL} we formulate the main results. The representation formulas for the solution of the initial value problem \eqref{01}, \eqref{02} are obtained in the general case: for arbitrary multi-order $\mathcal{B}$ and matrix-valued operator $\mathcal{F}(A).$ The main idea of the method used to obtain the representation formula is demonstrated for clarity first in the case $m=2,$ and then for arbitrary $m \ge 2.$ Note that some particular representation formulas were obtained in \cite{Umarov2024} in the case of systems of pseudo-differential operators. These results are also extended to differential-operator case.  Finally, in Section \ref{sec_appl} we discuss some applications and examples.

\section{Preliminaries}
\label{sec_prel}

\subsection{Fractional derivatives}
By definition the Riemann-Liouville fractional derivative of order $\beta \in (0,1)$ of a function defined on $[0,\infty),$ is the integral
\begin{equation} \label{RL_1}
D^{\beta}_{+} f(t) = \frac{1}{\Gamma(1-\beta)} \frac{d}{dt}\int_0^t \frac{f(\tau)d\tau}{(t-\tau)^{\beta}}.
\end{equation}
subject to exist, where $\Gamma(s)$ is the Euler's gamma-function.
Similarly, if $0 < \beta<1,$ then the Caputo derivative is defined by the integral
\begin{equation} \label{CD_1}
_a D^{\beta}_{\ast} f(t) = \frac{1}{\Gamma(1-\beta)} \int_0^t \frac{f^{\prime}(\tau)d\tau}{(t-\tau)^{\beta}}.
\end{equation}
subject to exist.

\noindent
The Laplace transform of the Riemann-Liouville and Caputo derivatives are 
\begin{align}
\label{RL_LT}
L[D_+^{\beta} f](s) & = s^{\beta}L[{f}](s) - (J^{1-\beta} f)(0), \\ 
L[D_{\ast}^{\beta} f](s) &= s^{\beta} L[f](s) -
f(0) s^{\beta-1},
\label{CD_LT}
\end{align}
respectively. We will use these formulas in the vector form. 
Namely, for a vector-valued function $F(t) =(f_1(t),\dots,f_m(t)),$ we have
\begin{align}
\label{RL_LT_01}
L[D_+^{\mathcal{B}} F(t)](s) & = s^{\mathcal{B}} L[F(t)](s) - (\mathcal{J}^{1-\mathcal{B}} F) (0) \notag \\
&= \Big(s^{\beta_1}L[{f_1}](s) - (J^{1-\beta_1} f_1)(0), \dots,  \notag \\ 
& s^{\beta_m}L[{f_m}](s) 
      - (J^{1-\beta_m} f_m)(0)\Big),
     \end{align}
where 
\begin{equation}
\label{FrInt_1}
(\mathcal{J}^{1-\mathcal{B}} \mathcal{F})(t) = \Big((J^{1-\beta_1}f_1)(t), \dots , (J^{1-\beta_m}f_m)(t) \Big).
\end{equation}
with fractional integrals 
\[
(J^{1-\beta_j} f)(t) = \frac{1}{\Gamma(1-\beta_j)}\int_0^t (t-\tau)^{\beta_j} f(\tau) d \tau, \quad j=1,\dots, m.
\]     
Similarly,
     \begin{align}
L[D_{\ast}^{\mathcal{B}} F(t)](s) &=  s^{\mathcal{B}} L[F(t)](s) -  F (0) s^{\mathcal{B}-1}  \notag \\ 
&=  \Big( s^{\beta_1} L[f_1](s) - f_1(0) s^{\beta_1-1},\dots,  \notag \\ 
& s^{\beta_m} L[f_m](s) 
       - f_m(0) s^{\beta_m-1}\Big). \label{CD_LT_01}
\end{align}
In these formulas 
\[
L[\mathcal{D}^{\mathcal{B}}F(t)](s)= \Big( L[\mathcal{D}^{\beta_1}f_1](s), \dots,  L[\mathcal{D}^{\beta_m} f_m](s) \Big)
\]
 for both operators $\mathcal{D}=D_+$ and $\mathcal{D}=D_{\ast}.$

\subsection{Matrix-valued functions}
Let $Z$ be a square matrix of size $m$ with the Jordan normal form  
\begin{equation}
\label{F1}
\mathbb{J}= M^{-1}Z M=\Lambda+N,
\end{equation}
 where $M$ is an invertible transformation matrix, $\Lambda$ is a diagonal matrix with eigenvalues on the diagonal and $N$ is the nilpotent matrix. Suppose $J_{\ell}, \ \ell = 1,\dots, L,$ are Jordan blocks of $Z$ 
and $\|Z\|$ is the matrix norm of $Z.$ Then for a function 
$g(z),$ analytic in a neighborhood of $|z| \le \|Z\|$ one has the spectral representation
\begin{align*}
g(\mathbb{J}) &= \sum_{k=0}^{m-1} \frac{N^k}{k!}  \int\limits_{\sigma(Z)} g^{(k)}(\lambda) dE_{\lambda}\\
&= \sum_{\ell=1}^L \sum_{k=0}^{m_{\ell}} \frac{N^k}{k!} g^{(k)}(\lambda_{\ell})\mathbb{P}_{\lambda_{\ell}},
\end{align*}
where the spectral measure $dE_{\lambda}$ is formed by projection operators $\mathbb{P}_{\lambda_{\ell}},$ 
determined by eigenvalues $\lambda_{\ell}, \ell=1,\dots,L,$ of the matrix $Z$ of multiplicity $m_{\ell}, \ell=1,\dots,L.$ In the explicit form this means
\begin{equation}
\label{m01}
g(\mathbb{J}) =
\begin{bmatrix}
g(J_1) & 0 & \dots & 0\\
0.    & g(J_2) & \dots & 0\\
\dots & \dots & \dots & \dots\\
0.    &  0.  &.  \dots    &  g(J_L)
\end{bmatrix},
\end{equation}
where
\begin{equation}
\label{m02}
g(J_{\ell}) =
\begin{bmatrix}
g(\lambda_{\ell}) & \frac{g^{\prime}(\lambda_{\ell})}{1!} &  \frac{g^{\prime \prime}(\lambda_{\ell})}{2!} & \dots &  \frac{g^{(m_{\ell}-2)}(\lambda_{\ell})}{(m_{\ell}-2)!}   & \frac{g^{(m_{\ell}-1)}(\lambda_{\ell})}{(m_{\ell}-1)!}  \\
0 & g(\lambda_{\ell}) & \frac{g^{\prime}(\lambda_{\ell})}{1!} &   \dots & \frac{g^{(m_{\ell}-3)}(\lambda_{\ell})}{(m_{\ell}-3)!}  & \frac{g^{(m_{\ell}-2)}(\lambda_{\ell})}{(m_{\ell}-2)!} \\
\dots & \dots & \dots & \dots & \dots & \dots \\
\dots & \dots & \dots & \dots & \dots & \dots \\
0      &    0    &     0    & \dots &      g(\lambda_{\ell})    &\frac{g^{\prime}(\lambda_{\ell})}{1!}\\
0.     &.   0.   &.    0.   & \dots &             0                       & g(\lambda_{\ell})
\end{bmatrix}, \quad \ell=1 \dots, L.
\end{equation}
It follows from \eqref{F1} and \eqref{m01} that
\[
g(Z) = M g(\mathbb{J}) M^{-1}.
\]


\subsection{Classical ML functions}
The classical two-parameter  Mittag-Leffler function is defined by
\[
E_{\beta, \nu}(z)= \sum_{n=0}^{\infty} \frac{z^n}{\Gamma(\beta n +\nu)}, \quad z \in \mathbb{C},
\]
where $\mathbb{C}$ is the set of complex numbers and parameters $\beta >0, \nu >0.$
This function plays an important role in the theory of fractional order differential equations. For various properties of the Mittag-Leffler function we refer the reader to sources \cite{GKMR,ML} and references therein. Here we only mention some properties of $E_{\beta, \nu}(z)$ used in the current paper. 
The function $E_{\beta, \nu}(z)$ is an entire function of order $1/\beta$ and recovers the exponential function $\exp(z)$ when $\beta =\nu=1.$ 
It is known \cite{GKMR,ML} that for $0 < \beta < 2, \ \nu \in \mathbb{C}, $ the Mittag-Leffler function $E_{\beta, \nu} (z)$
has asymptotic behavior 
\begin{equation}
\label{asymp01}
E_{\beta, \nu}(z)  \sim \frac{1}{\beta} z^{(1-\nu)/\beta}\exp(z^{1/\beta} ), \ |z| \to \infty, \quad \text{if } \ \  \frac{\beta \pi}{2} < |\arg(z)| < \min\{\pi, \beta \pi\},
\end{equation}
and 
\begin{equation}
\label{asymp2}
E_{\beta, \nu} (z) \sim 1/|z|,  |z| \to \infty, \quad  \text{if } \ \ \min\{\pi, \beta \pi\} < | \arg(z)|  \le \pi.
\end{equation}
For derivatives of $E_{\beta, \nu} (z^{\beta})$ the following formulas are valid:
\begin{equation}
\frac{d^k}{dz^k}\Big[ z^{\nu-1} E_{\beta, \nu} (z^{\beta}) \Big]  = z^{\nu - k -1} E_{\beta, \nu-k} (z^{\beta}), \ Re(\nu) > k, \ k = 1,2,\dots.
\end{equation}
Consider the function $E_{\beta, \nu}(\mu t^{\beta}),$ with a parameter $\mu \in \mathbb{C}.$ This function plays an important role in the theory of fractional order differential equations. For the Laplace transform of this function and its derivatives the following formulas hold \cite{SKM93,Umarov_book2015,GKMR}:
\begin{equation}
\label{L1}
L[t^{\nu-1}E_{\beta, \nu}(\mu t^{\beta})](s) = \frac{s^{\beta-\nu}}{s^{\beta}-\mu}, \quad s> [Re(\mu)]^{1/\beta},
\end{equation}
\begin{equation}
\label{L3}
L\left[ \frac{t^{k\beta +\nu-1}}{k!}E^{(k)}_{\beta, \nu}(\mu t^{\beta}) \right](s) = \frac{s^{\beta-\nu}}{(s^{\beta}-\mu)^{k+1}}, \quad s> [Re(\mu)]^{1/\beta}, k=1,2,\dots,
\end{equation}
where $E_{\beta, \nu}^{(k)}(z) = \frac{d^{k}}{dz^k}E_{\beta, \nu}(z).$ 
In particular, if $\nu=1,$ then one gets
\begin{equation}
\label{L2}
L\left[ \frac{t^{k\beta}}{k!}E^{(k)}_{\beta}(\mu t^{\beta}) \right](s) = \frac{s^{\beta-1}}{(s^{\beta}-\mu)^{k+1}}, \quad s> [Re(\mu)]^{1/\beta}, k=0,1,\dots,
\end{equation}
and if $\nu=\beta$ in \eqref{L3},
\begin{equation}
\label{L4}
L\left[ \frac{t^{k\beta +\beta-1}}{k!}E^{(k)}_{\beta, \beta}(\mu t^{\beta}) \right](s) = \frac{1}{(s^{\beta}-\mu)^{k+1}}, \quad s> [Re(\mu)]^{1/\beta}, k=0,1,\dots,\end{equation}

The convolution of functions $f(t), \ g(t), \ t \ge 0,$ is defined by
\begin{equation}
\label{conv}
(f\ast g)(t) = \int_0^t f(\tau) g(t-\tau) d\tau.
\end{equation}
The following lemmas will be used in our further analysis. 

\begin{lem} \label{lemma_1}
For  $0<\beta_1<\beta_2, \ \nu>0,$ and $k = 0,1,\dots,$ the following relations hold:

\begin{enumerate}
\item[(a)]
\begin{equation}
\label{rel_new}
\Big( I-\mu J^{\beta} \Big)^k \Big[ \frac{t^{k \beta +\nu-1}}{k!} E^{(k)}_{\beta, \nu}(  \mu t^{\beta} ) \Big] = J^{k \beta} E_{\beta, \nu}( \mu t^{\beta} ),
\end{equation}
\item[(b)]
\begin{equation}
\label{rel_new+}
\Big( J^{\beta_2-\beta_1}-\mu J^{\beta_2} \Big)^k \Big[ \frac{t^{k \beta_1 +\nu-1}}{k!} E^{(k)}_{\beta_1, \nu}(  \mu t^{\beta_1} ) \Big] = J^{k \beta_2} E_{\beta_1, \nu}( \mu t^{\beta_1} ),
\end{equation}
\end{enumerate}
where $I$ is the identity operator and $J^{\beta}$ is the fractional integral of order $\beta.$
\end{lem}
{\it Proof}  (a) To proof this statement we show that the Laplace transforms of both sides in \eqref{rel_new} coincide. 
Indeed, applying the Laplace transform to the left side of \eqref{rel_new}, we have
\[
\left( 1-\frac{\mu}{s^{\beta}} \right)^k \frac{s^{\beta-\nu}}{(s^{\beta}-\mu)^{k+1}} = \frac{s^{\beta-\nu}}{s^{k \beta}(s^{\beta}-\mu)}, \quad s> [Re(\mu)]^{1/\beta}.
\]
This is obviously the Laplace transform of the right hand side of \eqref{rel_new}, as well. 

(b) Similarly, the Laplace transform of the left hand side of \eqref{rel_new+} is

\[
\left( \frac{1}{s^{\beta_2-\beta_2}} - \frac{\mu}{s^{\beta_2}} \right)^k \frac{s^{\beta_1-\nu}}{(s^{\beta_1}-\mu)^{k+1}} = \frac{s^{\beta_1-\nu}}{s^{k\beta_2}(s^{\beta_1}-\mu)}, \quad s> [Re(\mu)]^{1/\beta_1},
\]
which is the Laplace transform of the right hand side of \eqref{rel_new+}, as well.
\eproof

\begin{lem} \label{lemma_2}
For any $\beta_1>0, \beta_2>0, \nu>0,$ and parameters $\mu_1, \mu_2 \in \mathbb{C}$ the following relations hold:
\begin{enumerate}
\item[(i)]
\begin{equation}
\label{rel_new2}
(I-\mu_2 J^{\beta_2})^{-1} \Big[t^{\nu-1}E_{\beta_1,\nu}(\mu_1 t^{\beta_1}) -\mu_1 J^{\nu} (t^{\beta_1-1} E_{\beta_1,\beta_1}(\mu_1 t^{\beta_1})\Big] = t^{\nu-1}E_{\beta_2, \nu}(\mu_2 t^{\beta_2}),
\end{equation}
\item[(ii)]
\begin{equation}
\label{rel_new3}
(I-\mu_2 J^{\beta_2})^{-1} J^{\beta_2}[t^{\nu-1}E_{\beta_1, \nu}(\mu_1 t^{\beta_1})] = \Big(t^{\nu-1}E_{\beta_1,\nu}(\mu_1 t^{\beta_1}) \Big) \ast \Big( t^{\beta_2-1}E_{\beta_2, \beta_2}(\mu_2 t^{\beta_2}) \Big),
\end{equation}
\end{enumerate}
where $"\ast``$ is the convolution operation.
\end{lem}
{\it Proof} (i)
Again we show that the Laplace transforms of both sides in equations \eqref{rel_new2} and \eqref{rel_new3} coincide. For the Laplace transform of the left hand side of \eqref{rel_new2} we have
\begin{align*}
\mathcal{L} & \Big[ (I-\mu_2 J^{\beta_2})^{-1} [t^{\nu-1}E_{\beta_1, \nu}(\mu_1 t^{\beta_1}) -\mu_1 J^{\nu} E_{\beta_1,\beta_1}(\mu_1 t^{\beta_1})] \Big]
\\
&= (1-\frac{\mu_2}{s^{\beta_2}})^{-1} \Big( \frac{s^{\beta_1-\nu}}{s^{\beta_1}-\mu_1} - \frac{\mu_1}{s^{\nu}(s^{\beta_1}-\mu_1)} \Big)
\\
&= \frac{s^{\beta_2}}{s^{\beta_2}-\mu_2} \ \frac{s^{\beta_1} - \mu_1}{s^{\nu}(s^{\beta_1}-\mu_1)} =\frac{s^{\beta_2-\nu}}{s^{\beta_2}-\mu_2},
\end{align*}
where
\[
 s> \max \Big\{ [Re(\mu_1)]^{1/\beta_1},  [Re(\mu_2)]^{1/\beta_2} \Big\}.
 \]
This is obviously the Laplace transform of the right hand side of \eqref{rel_new2}, as well. 

(ii) Similarly, the Laplace transform of the left hand side of \eqref{rel_new3} is
\begin{align*}
\mathcal{L}  & \Big[ (I-\mu_2 J^{\beta_2})^{-1} J^{\beta_2}[t^{\nu-1}E_{\beta_1, \nu}(\mu_1 t^{\beta_1})] \Big] 
= (1-\frac{\mu_2}{s^{\beta_2}})^{-1}  \frac{s^{\beta_1-\nu}}{s^{\beta_2}(s^{\beta_1}-\mu_1)}
\\
&=  \frac{s^{\beta_2}}{s^{\beta_2}-\mu_2} \frac{s^{\beta_1-\nu}}{s^{\beta_2}(s^{\beta_1}-\mu_1)} = \frac{s^{\beta_1-\nu}}{s^{\beta_1}-\mu_1} \ \frac{1}{s^{\beta_2}-\mu_2}.
\end{align*}
On the other hand, due to the convolution theorem,  the Laplace transform of the right hand side of \eqref{rel_new3} also results in the same expression.
\eproof

\subsection{Matrix-valued ML functions}

Since $E_{\beta, \nu}(z)$ is an entire function, in accordance with \eqref{m01} and \eqref{m02}, for a matrix $Z$ one can introduce a matrix-valued version of the Mittag-Leffler function as
\begin{equation}
\label{ML}
\mathbb{E}_{\beta, \nu}( Z) = M \mathbb{E}_{\beta, \nu}( \mathbb{J}) M^{-1} =M \mathbb{E}_{\beta, \nu} ( {\Lambda} +N ) M^{-1},
\end{equation}
where
\begin{equation}
\label{m1}
\mathbb{E}_{\beta, \nu}( \mathbb{J}) =
\begin{bmatrix}
\mathbb{E}_{\beta, \nu}(J_1) & 0 & \dots & 0\\
0.    & \mathbb{E}_{\beta, \nu}(J_2) & \dots & 0\\
\dots & \dots & \dots & \dots\\
0.    &  0.  &.  \dots    &  \mathbb{E}_{\beta, \nu}(J_L)
\end{bmatrix},
\end{equation}
with the block matrices 
\begin{align}
\mathbb{E}_{\beta, \nu}( J_{\ell}) &=
\begin{bmatrix}
E_{\beta, \nu}(\lambda_{\ell}) & \frac{E_{\beta, \nu}^{\prime}(\lambda_{\ell})}{1!}  & \dots &  \frac{E_{\beta, \nu}^{(m_{\ell}-2)}(\lambda_{\ell})}{(m_{\ell}-2)!}   & \frac{E_{\beta, \nu}^{(m_{\ell}-1)}(\lambda_{\ell})}{(m_{\ell}-1)!}  \\
0 & E_{\beta, \nu}(\lambda_{\ell})  &   \dots & \frac{E_{\beta, \nu}^{(m_{\ell}-3)}(\lambda_{\ell})}{(m_{\ell}-3)!}  & \frac{E_{\beta, \nu}^{(m_{\ell}-2)}(\lambda_{\ell})}{(m_{\ell}-2)!} \\
\dots & \dots  & \dots & \dots & \dots \\
\dots & \dots & \dots & \dots & \dots \\
0      &    0       & \dots &      E_{\beta, \nu}(\lambda_{\ell})    &\frac{E_{\beta, \nu}^{\prime}(\lambda_{\ell})}{1!}\\
0     &   0      & \dots &             0                       & E_{\beta, \nu}(\lambda_{\ell})
\end{bmatrix},
\label{m03}
\end{align}
corresponding to (algebraic) eigenvalues  $\lambda_{\ell}, \ell=1 \dots, L,$ of the matrix $Z.$ It is not hard to verify that using formulas \eqref{L1}, \eqref{L3}, and \eqref{m03} one obtains the Laplace transforms of  the matrix-valued function $ t^{ \nu -1}\mathbb{E}_{\beta, \nu}(t^{\beta}J_{\ell}):$       
\begin{equation}
\label{m4}
L[t^{\nu-1}\mathbb{E}_{\beta, \nu}(t^{\beta}J_{\ell})](s)=
\begin{bmatrix}
\frac{s^{\beta-\nu}}{s^{\beta}-\lambda_{\ell}} & \frac{s^{\beta-\nu}}{(s^{\beta}-\lambda_{\ell})^2} & \dots & \frac{s^{\beta-\nu}}{(s^{\beta}-\lambda_{\ell})^{m_{\ell}}} \\ 
0  & \frac{s^{\beta-\nu}}{s^{\beta}-\lambda_{\ell}} & \dots & \frac{s^{\beta-\nu}}{(s^{\beta}-\lambda_{\ell})^{m_{\ell}-1}} \\
\dots & \dots & \dots & \dots \\
0 & 0 & \dots & \frac{s^{\beta-\nu}}{s^{\beta}-\lambda_{\ell}}
\end{bmatrix}, \quad \ell=1,\dots L.
\end{equation}

\subsection{Vector-indiced matrix-valued ML functions}

Let $\mathcal{B}=(\beta_1, \dots, \beta_m)$ and ${\mathcal{V}}=(\nu_1,\dots,\nu_m)$ be vector indices with components $\beta_j >0, \ \nu_j >0, \ j=1,\dots,m. $ For a diagonal matrix $D$ with diagonal entries $d_1,\dots,d_m,$ we use the notation
\[
D= \mbox{diag} (d_1,\dots, d_m).
\]
\begin{definition}
\label{def_MVML} Let $Z$ be a square matrix of size $m\times m$ with complex entries. A vector-indiced matrix-valued Mittag-Leffler function denoted by $E_{\mathcal{B}, \mathcal{V}}(z),$ is defined by
\begin{equation}
\label{ML_1}
E_{\mathcal{B}, \mathcal{V}}(Z) = \sum_{n=0}^{\infty} \Big(I \Gamma(n \mathcal{B}+\mathcal{V}) \Big)^{-1} Z^n,
\end{equation}
where $I \Gamma(n \mathcal{B}+\mathcal{V}) =\text{diag} \Big(\Gamma(n\beta_1 + \nu_1), \dots, n \Gamma(\beta_m+\nu_m)\Big).$
\end{definition}

The vector-indiced matrix-valued ML function $E_{\mathcal{B}, \mathcal{V}}(Z)$ generalizes the classical and above considered matrix-valued ML functions. Here are some examples:
\begin{enumerate}
\item Let $m=1$ and $\beta_1=\beta, \nu_1=\nu.$ Then we obtain the classical two-parameter ML function $E_{\beta, \nu}(z).$
\item Let $\mathcal{B}=(\beta,\dots,\beta)$ and $\mathcal{V}=(\nu_1,\dots,\nu).$ Then we obtain the matrix-valued ML function \cite{Garrapa,GKMR}
\begin{equation}
\label{MLmv}
E_{\beta, \nu} (Z) = \sum_{n=0}^{\infty} \frac{Z^n}{\Gamma (n \beta + \nu)}.
\end{equation}
Indeed, in this case
\begin{align*}
 \Big(I \Gamma(n \mathcal{B}+\mathcal{V}) \Big)^{-1} Z^n &= \Big[ \text{diag} (\Gamma(n\beta + \nu), \dots, \Gamma(n \beta+\nu))\Big]^{-1} Z^n 
 \\
 &=  \text{diag} \left(\frac{1}{\Gamma(n\beta + \nu)}, \dots, \frac{1}{\Gamma(n \beta+\nu)} \right) Z^n 
 \\
 &=\frac{1}{\Gamma(n\beta + \nu)} Z^n. 
  \end{align*}
Therefore, in this case \eqref{ML_1} reduces to \eqref{MLmv}.
\end{enumerate}

Let $\lambda_j, j=1,\dots, L,$ be eigenvalues of (algebraic) multiplicity $m_j$  of the matrix $Z=M(\Lambda+N)M^{-1},$ and let $J_{\ell}$ be the Jordan block of the Jordan canonical form $\Lambda + N$ corresponding to $\lambda_j.$ Then, it is not hard to see that
\begin{equation}
\label{m10}
\mathbb{E}_{\mathcal{B}, \mathcal{V}}({\Lambda+N}) =
\begin{bmatrix}
\mathbb{E}_{\mathcal{B}_1, \mathcal{V}_1}(J_1) & 0 & \dots & 0\\
0    & \mathbb{E}_{\mathcal{B}_2, \mathcal{V}_2}(J_2) & \dots & 0\\
\dots & \dots & \dots & \dots\\
0    &  0  &  \dots    &  \mathbb{E}_{\mathcal{B}_L, \mathcal{V}_L}(J_L)
\end{bmatrix},
\end{equation}
where $\mathcal{B}_{\ell} =(\beta_{M_{\ell}+1}, \dots, \beta_{M_{\ell}+m_{\ell}}), \ \mathcal{V}_{\ell} =(\nu_{M_{\ell}+1}, \dots, \nu_{M_{\ell}+m_{\ell}}),$ and \\$M_{\ell} = m_1 + \dots + m_{\ell-1},$ with blocks
\begin{align}
\mathbb{E}_{\mathcal{B}_{\ell}, \mathcal{V}_{\ell}}( J_{\ell}) 
&=
\begin{bmatrix}
E_{\beta_{M_{\ell}}+1, \nu_{M_{\ell}}+1}(\lambda_{\ell}) & \frac{E_{\beta_{M_{\ell}}+1, \nu_{M_{\ell}}+1}^{\prime}(\lambda_{\ell})}{1!}  & \dots  
& \frac{E_{\beta_{M_{\ell}}+1, \nu_{M_{\ell}}+1}^{(m_{\ell}-1)}(\lambda_{\ell})}{(m_{\ell}-1)!}  \\
0 & E_{\beta_{M_{\ell}}+2, \nu_{M_{\ell}}+2}(\lambda_{\ell})  &   \dots 
& \frac{E_{\beta_{M_{\ell}}+2, \nu_{M_{\ell}}+2}^{(m_{\ell}-2)}(\lambda_{\ell})}{(m_{\ell}-2)!} \\
\dots & \dots  & \dots & \dots & \\ 
0     &   0      & \dots                      
& E_{\beta_{M_{\ell}}+m_{\ell}, \nu_{M_{\ell}}+m_{\ell}}(\lambda_{\ell})
\end{bmatrix}.
\label{ml03}
\end{align}

Now suppose that 
\[
\mathcal{B} = (\underset{\mbox{$m_1$ times}}{\underbrace{ \beta_1, \dots, \beta_1}}, \underset{\mbox{$m_2$ times}}{\underbrace{ \beta_2, \dots, \beta_2}}, \dots, \underset{\mbox{$m_L$ times}}{\underbrace{ \beta_L, \dots, \beta_L}}),
\]
and
\[
\mathcal{V} = (\underset{\mbox{$m_1$ times}}{\underbrace{ \nu_1, \dots, \nu_1}}, \underset{\mbox{$m_2$ times}}{\underbrace{ \nu_2, \dots, \nu_2}}, \dots, \underset{\mbox{$m_L$ times}}{\underbrace{ \nu_L, \dots, \nu_L}}).
\]
Then, it follows from \eqref{ml03} that for each $\ell = 1, \dots, L:$
\begin{align}
\notag
&\mathbb{E}_{\mathcal{B}_{\ell}, \mathcal{V}_{\ell}}(I t^{\mathcal{B}_{\ell}} J_{\ell}) \\
\hspace{-1cm}
&=
\begin{bmatrix}
E_{\beta_{\ell}, \nu_{\ell}}(\lambda_{\ell}t^{\beta_{\ell}}) & \frac{t^{\beta_{\ell}}E_{\beta_{\ell}, \nu_{\ell}}^{\prime}(\lambda_{\ell}t^{\beta_{\ell}})}{1!}  & \dots &  \frac{t^{(m_{\ell}-2){\beta_{\ell}}}E_{\beta_{\ell}, \nu_{\ell}}^{(m_{\ell}-2)}(\lambda_{\ell}t^{\beta_{\ell}})}{(m_{\ell}-2)!}   & \frac{t^{(m_{\ell}-1)\beta_{\ell}}E_{\beta_{\ell}, \nu_{\ell}}^{(m_{\ell}-1)}(\lambda_{\ell}t^{\beta}_{\ell})}{(m_{\ell}-1)!}  \\
0 & E_{\beta_{\ell}, \nu_{\ell}}(\lambda_{\ell}t^{\beta_{\ell}})  &   \dots & \frac{t^{(m_{\ell}-3)\beta}E_{\beta_{\ell}, \nu_{\ell}}^{(m_{\ell}-3)}(\lambda_{\ell}t^{\beta_{\ell}})}{(m_{\ell}-3)!}  & \frac{t^{(m_{\ell}-2)\beta}E_{\beta_{\ell}, \nu_{\ell}}^{(m_{\ell}-2)}(\lambda_{\ell}t^{\beta_{\ell}})}{(m_{\ell}-2)!} \\
\dots & \dots  & \dots & \dots & \dots \\
\dots & \dots & \dots & \dots & \dots \\
0      &    0       & \dots &      E_{\beta_{\ell}, \nu_{\ell}}(\lambda_{\ell}t^{\beta_{\ell}})    &\frac{t^{\beta_{\ell}}E_{\beta_{\ell}, \nu_{\ell}}^{\prime}(\lambda_{\ell}t^{\beta_{\ell}})}{1!}\\
0     &   0      & \dots &             0                       & E_{\beta_{\ell}, \nu_{\ell}}(\lambda_{\ell}t^{\beta_{\ell}})
\end{bmatrix},
\label{ml300}
\end{align}
where $I t^{\mathcal{B}_{\ell}} = \mbox{diag} (t^{\beta_{\ell}}, \dots, t^{\beta_{\ell}})$ is a diagonal matrix of size $m_{\ell} \times m_{\ell}.$ 

The latter implies

\begin{equation}
\label{m400}
L[It^{\mathcal{B}_{\ell}-\mathcal{V}_{\ell}}\mathbb{E}_{\mathcal{B}_{\ell}, \mathcal{V}_{\ell}}(t^{\mathcal{B}_{\ell}}J_{\ell})](s)=
\begin{bmatrix}
\frac{s^{\beta_{\ell}-\nu_{\ell}}}{s^{\beta_{\ell}}-\lambda_{\ell}} & \frac{s^{\beta_{\ell}-\nu_{\ell}}}{(s^{\beta_{\ell}}-\lambda_{\ell})^2} & \dots & \frac{s^{\beta_{\ell}-\nu_{\ell}}}{(s^{\beta_{\ell}}-\lambda_{\ell})^{m_{\ell}}} \\ 
0  & \frac{s^{\beta_{\ell}-\nu_{\ell}}}{s^{\beta_{\ell}}-\lambda_{\ell}} & \dots & \frac{s^{\beta_{\ell}-\nu_{\ell}}}{(s^{\beta_{\ell}}-\lambda_{\ell})^{m_{\ell}-1}} \\
\dots & \dots & \dots & \dots \\
0 & 0 & \dots & \frac{s^{\beta-1}}{s^{\beta_{\ell}}-\lambda_{\ell}}
\end{bmatrix}, \quad \ell=1,\dots L.
\end{equation}

We note that, in general,  $M \mathbb{E}_{\mathcal{B}, \mathcal{V}}({\Lambda+N}) M^{-1}$ is not the same as ${E}_{\mathcal{B}, \mathcal{V}}({Z})$ unless vector-indices $\mathcal{B}$ and $\mathcal{V}$ have equal components. Indeed, using the equality 
\[
Z^n = M (\Lambda+N)^n M^{-1},
\] 
one obtains
\[
{E}_{\mathcal{B}, \mathcal{V}}({Z}) = \sum_{n=0}^{\infty} \Big(I \Gamma(n \mathcal{B}+\mathcal{V}) \Big)^{-1} \ M \ (\Lambda+N)^n \ M^{-1} \neq M \mathbb{E}_{\mathcal{B}, \mathcal{V}}({\Lambda+N}) M^{-1},
\]
since matrices $\Big(I \Gamma(n \mathcal{B}+\mathcal{V}) \Big)^{-1} $ and $M$ do not commute. It is not hard to verify that these matrices commute if and only if vectors $\mathcal{B}$ and $\mathcal{V}$ have equal components. In this case the following theorem holds.

\begin{thm} 
\label{thm_gen}
Let $\mathcal{B}=(\beta, \dots, \beta), \ 0<\beta \le1,$ and $\mathcal{V}=(\nu, \dots, \nu), \nu>0.$ Then
\begin{enumerate}
\item for the matrix-valued ML function ${E}_{\mathcal{B}, \mathcal{V}}({Z})$ the following representation is valid:
\[
{E}_{\mathcal{B}, \mathcal{V}}({Z}) = M {E}_{\mathcal{B}, \mathcal{V}}({\Lambda+N}) M^{-1};
\]
\item
the following Laplace transform 
formula holds:
\begin{equation}
\label{LT}
\mathcal{L}[It^{\mathcal{V}- {\bf 1}} {E}_{\mathcal{B}, \mathcal{V}}({It^{\mathcal{B}}Z})](s) = Is^{\mathcal{B}-\mathcal{V}}(I s^{\mathcal{B}}-Z)^{-1},
\end{equation}
where $\mathcal{V}-{\bf 1}=(\nu-1,\dots,\nu-1)$ and $It^{\mathcal{B}}=\mbox{diag}(t^{\beta}, \dots, t^{\beta}).$
\end{enumerate}
\end{thm}

{\it Proof}  \ We need to prove only part 2. Using the definition \eqref{ML_1} of the ML function, we have

\begin{equation}
\label{ML_71}
\mathcal{L}\Big[It^{\mathcal{V}-{\bf 1}}E_{\mathcal{B}, \mathcal{V}}(I t^{\mathcal{B}}Z) \Big](s)  = \sum_{n=0}^{\infty} \Big( I \Gamma(n \mathcal{B}+\mathcal{V}) \Big)^{-1} \mathcal{L}\Big[ I t^{n \mathcal{B}+ \mathcal{V} - {\bf 1}} \Big] Z^n,
\end{equation}
since matrices $It^{\mathcal{V}-{\bf 1}}$  and $I \Gamma(n \mathcal{B}+\mathcal{V}) \Big)^{-1}$ commute under the conditions to $\mathcal{B}$ and $\mathcal{V}.$ Further, using the well known relation $\mathcal{L}[t^{\rho}](s)= \Gamma(\rho+1)/s^{\rho+1},$ we obtain
\begin{equation}
\label{ML_72}
\mathcal{L}\Big[It^{\mathcal{V}-{\bf 1}}E_{\mathcal{B}, \mathcal{V}}(I t^{\mathcal{B}}Z) \Big](s)  = Is^{\mathcal{-V}}\sum_{n=0}^{\infty}\Big[ I s^{- n \mathcal{B}} Z^n \Big] =  Is^{\mathcal{B}-\mathcal{V}}(I s^{\mathcal{B}}-Z)^{-1},
\end{equation}
completing the proof.
\eproof

\subsection{Matrix-valued operators with singular symbols}
\label{fcalculus}

In this section we describe matrix-valued operators $\mathcal{F}(A)$ on the right hand side of system \eqref{01}. 
Let $A$ be a closed linear operator with a domain ${\cal{D}}\,
(A)$ dense in a reflexive Banach space ${X}$ and a nonempty spectrum $\sigma (A) \subset \mathbb{C}.$
Assume the entry $f_{jk}(A)$ of the matrix-valued operator $\mathcal{F}(A)$ has the symbol  $f_{jk}(z), \ z \in \mathbb{C},$ analytic in an open connected domain $G \subset \mathbb{C}.$ If $\sigma(A)$ is bounded and $G$
contains $\sigma(A)$ then one can define the operator $f_{jk}(A)$ as (see e.g. \cite{DalKrein,DSh})
\begin{equation}\label{1}
f_{jk}(A) = \int_{\gamma} f_{jk}(\zeta) {\mathcal{R}} (\zeta, A)  d \zeta, \quad j,k=1,\dots,m,
\end{equation}
where $\gamma$ is a contour in $G$ containing $\sigma(A),$ and
$\mathcal{R}(\zeta, A), ~ \zeta \in \mathbb{C} \setminus \sigma(A),$ is
the resolvent operator of $A.$ Representation \eqref{1} is not valid if  $f_{jk}(z)$ has singularities on the spectrum $\sigma(A).$

In the case when $f$ has singular points in the spectrum $\sigma(A)$ of the operator $A$ the corresponding operator $f(A)$ 
can be constructed as follows. Denote by $\sing(f)$ the set of singular points of $f$ on $\sigma(A).$ Let $\mathbb{D}$ be an open
set in $\mathbb{C}$ containing $\sigma(A).$ In particular, if $\sigma(A) = \mathbb{C}$ then $\mathbb{D} = \mathbb{C}$ as well. Consider an open set $G \subset \mathbb{D} \setminus \sing(f).$  Let $0 < r \leq +\infty$ and $\mu < r$. Denote by
$X_{\mu}$ the set of elements $x \in \cap_{k\ge 1}{\mathcal{D}}(A^k)$ satisfying the inequalities
$\|A^{k} x \| \leq C \mu^{k}\|x\|$ for all $k = 1, 2, ...,$ with a constant $C>0$ not depending on $k$. 
A sequence of elements $x_{n} \in X_{\mu},\,\,n=1, 2, ...,$ is said to converge to an element $x_{0} \in X_{ \mu}$ if 
$\|x_{n} - x_{0}\| \rightarrow 0,$ as $n \rightarrow \infty$. It is easy to see that $X_{\mu_1} \subset X_{\mu_2},$ if $\mu_1 < \mu_2,$ and this inclusion is continuous. 
Denote by $X_{A,r}$ the inductive limit of spaces $X_{\mu}$ as $\mu \to r,$ i.e.
 \[
 X_{A, r} = \underset{\mu \to r}{ \mbox{ind-lim}} \ X_{\mu},
 \]
 meaning that $X_{A,r} = \cup_{0<\mu < r} X_{\mu}$ with the strongest topology.  
For basic notions of topological vector spaces including inductive and
 projective limits we refer the reader to \cite{Robertson}. 
 The space $X_{A, r}$ is called a space of exponential vectors of type $r$ (see e.g. \cite{Radino,Um98}) associated with the operator $A.$  
 
 Let $A_{\lambda} = A-\lambda I,$ where $\lambda \in G,$ and denote by $X_{A,\,G}$ the space whose elements are locally finite sums of  elements in $X_{A_{\lambda}, r}, \ r < \mbox{dist} (\lambda, \partial G),$
 with the corresponding topology. Here $\mbox{dist} (\lambda, \partial G)$ is the minimal distance between the point $\lambda$ and  boundary of the domain $G.$ By definition, any $u\in X_{A,\,G}$ has a
representation 
\[
u=\sum_{k=1}^{m_u} u_{\lambda_k}, \quad u_{\lambda_k} \in X_{A_{\lambda_k}, \, r},
\] 
where $\lambda_k \in G,$ and $m_u$ is a finite number. 

 Now we can define operators $f(A)$ with symbols $f(z)$ analytic in the domain $G.$ Recall that $f(z)$ may have singular points on the spectrum $\sigma(A),$ but $G$ does not contain singularities of $f(z).$ As an analytic function in $G,$ $f(z)$ 
 has the Taylor expansion 
 \[
 f_{\lambda}(z)=\sum_{n=0}^{\infty} \frac{f^{(n)}(\lambda)}{n!} (z-\lambda)^n, \quad \lambda \in G,
 \]
 convergent in any open disk $|z-\lambda| < r,$ where $r < \mbox{dist}(\lambda, \partial G).$ Therefore, 
 the operator $f_{\lambda}(A)$  defined as 
 \begin{equation}
 \label{f_lambda}
 f_{\lambda}(A)u_{\lambda}=\sum_{n=0}^{\infty} \frac{f^{(n)}(\lambda)}{n!} A_\lambda^n u_{\lambda}
  \end{equation}
on  elements $u_{\lambda} \in X_{A_{\lambda},\,r}$
 is well defined. Indeed, we have
 \begin{align}
 \| f_{\lambda}(A)u_{\lambda}\| &\le C \sum_{n=0}^{\infty} \frac{|f^{(n)}(\lambda)|}{n!} \|A_{\lambda}^n u_{\lambda} \| \notag \\
&\le  C \|u_{\lambda}\| \sum_{n=0}^{\infty} \frac{|f^{(n)}(\lambda)|}{n!}  \mu^n <\infty, \quad \mu<r.
\label{est}
  \end{align}
  Finally, for  an arbitrary $u \in X_{A,\,G}$ with the representation 
  \begin{equation}
  \label{repres}
  u = \sum_{\lambda \in G} u_{\lambda},\,\,\,u_{\lambda} \in X_{A_{\lambda},\,r},
  \end{equation}
 the operator $f(A)$ is defined by the formula
\begin{equation}
\label{calc1} 
f(A)u = \sum_{\lambda \in G} f_{\lambda}(A)u_{\lambda},
 \end{equation}
where $f_{\lambda}(A)u_{\lambda}$  is defined in \eqref{f_lambda}.
 Using estimate \eqref{est} and representation \eqref{calc1}, it is easy to show that the operator $f(A)$ is well defined on the space $X_{A,\,G}.$

Further, suppose that there exists a one-parameter family of
bounded invertible operators $E_{\lambda}: {X}\rightarrow {X}$ such
that
\begin{equation}
\label{ulambda} A_{\lambda} = E_{\lambda}A E^{-1}_{\lambda},
\,\,\lambda \in G.
\end{equation}
\begin{ex} 
\label{ex1}
Let $X=L_2 \equiv L_2(R)$ and $A = -i\frac{d}{dx}: L_2
\rightarrow L_2$  with domain
${\mathcal{D}} (A) = \{ v \in L_2: A v \in L_2 \}.$ Then for
operators $E_{\lambda}: v(x) \rightarrow e^{i\lambda x} v(x) $ we have
\begin{align*}
AE_{\lambda}v(x)&=-i \frac{d}{dx}(e^{i\lambda x}v(x)) = \lambda
e^{i\lambda x}v(x) - i e^{i\lambda x}\frac{dv}{dx}  \\
&= \lambda E_{\lambda} v(x) + E_{\lambda} A v(x),
\end{align*}
obtaining (\ref{ulambda}).
\end{ex} 
It follows from \eqref{ulambda} that
\[
A_{\lambda}^n = E_{\lambda} A^n E_{\lambda}^{-1},
\]
for all $n = 1,2,\dots,$ yielding  
\begin{equation}
\label{calc1.1} 
f(A)u = \sum_{\lambda \in G} \sum_{n=0}^{\infty}
  \frac{f^{(n)}(\lambda)}{n!} U_{\lambda}A^{n}U_{\lambda}^{-1}u_{\lambda}.
\end{equation}
Recall that here the sum with respect to $\lambda$ is finite.

The operator $f(A)$ defined in \eqref{calc1} maps $X_{A, G}$ to itself. Namely the mapping
\[
f(A): X_{A, G} \rightarrow X_{A, G}
\]
is continuous.
Indeed, let $u \in X_{A, G}$ has a representation $u=\sum_{\lambda}u_{\lambda}, ~
u_{\lambda} \in X_{A_{\lambda}, r}.$ Then
for $f(A)u$ we have the estimate
\begin{align}
\|A^{k}_{\lambda} f_{\lambda}(A) u_{\lambda} \|
&\le \sum_{n=0}^{\infty} \frac{|f^n(\lambda)|}{n!}\|(A-\lambda
I)^nA_{\lambda}^ku_{\lambda}\| 
\notag \\
&\leq \max_{|z-\lambda|<r} |f(z)|  \|A_{\lambda}^k u_{\lambda}\| \le C \mu^k \|u_{\lambda}\|,
\label{estimate1} 
\end{align}
 with some constant $C>0$ and $\mu <r.$ The latter means that $f_{\lambda}(A)u_{\lambda} \in X_{A_{\lambda}, r}.$ Therefore, $f(A)u$ has a representation 
$\sum_{\lambda}v_{\lambda},$ where $v_{\lambda}=f_{\lambda}(A)u_{\lambda} \in X_{A_{\lambda}, r},$ implying $f(A)u \in
X_{A, G}.$ The estimate (\ref{estimate1})
also implies the continuity of the mapping $f(A)$ in the topology of
$X_{A, G}.$
\begin{rem}
If the spectrum of the operator $A$ is discrete, then $X_{A_{\lambda}, r}$ consists of all linear combinations of eigenvectors and associated eigenvectors corresponding to eigenvalues $\lambda_k$ in the the disc $|\lambda - \lambda_k| < r,$ and the space $X_{A, G}$ is their locally finite sums.
\end{rem}

Finally, it follows from the construction above that a matrix-valued operator $\mathcal{F}(A)$ with the matrix-symbol $\mathcal{F}(z)=\{f_{kj}(z), \ k,j = 1, \dots, m \},$ analytic in the domain $G,$ is well defined on elements of the direct product space
\[
\mathcal{X}_{A,G} =X_{A, G} \underset{\mbox{$m$ times}}{\underbrace{\otimes \dots \otimes}} X_{A, G}, 
\]
with the corresponding direct product topology. Moreover, the mapping
\begin{equation}
\label{mapping}
\mathcal{F}(A): \mathcal{X}_{A,G} \to \mathcal{X}_{A,G}
\end{equation}
is continuous.

We note that the space $\mathcal{X}_{A,G}$ is relatively narrow. For example, if $A=(-i \frac{\partial}{\partial x_1}, \dots , -i \frac{\partial}{\partial x_n})$ acting in the space $L_2(\re^n),$ then the corresponding space $\mathcal{X}_{A,G}$ is the direct product of the space of functions analytic in $G \subset \re^n.$ However, the duality construction allows to expand the introduced spaces and consider wider classes of fractional order systems.
Let $X^{\ast}$ denote the dual of $X,$ and $A^{\ast}: X^{\ast} \to
X^{\ast}$ be the  operator adjoint to $A$. We denote by
$X^{'}_{A^{\ast}, G^{\ast}}$ the space
of linear continuous functionals defined on $X_{A,
G},$ with respect to weak convergence. In other words $X^{'}_{A^{\ast}, G^{\ast}}$ is the projective limit of spaces $X^{'}_{A^{\ast}_{\lambda}, r}$ which are dual to $X_{A_{\lambda}, r}$ with the coarsest topology. Continuing the example $A=(-i \frac{\partial}{\partial x_1}, \dots , -i \frac{\partial}{\partial x_n}),$ now one can see that the corresponding space $X^{'}_{A^{\ast}, G^{\ast}}$ gives rise to the space of analytic functionals (Sato's hyperfunctions; see e.g. \cite{Sato}).

For an analytic matrix-symbol $\mathcal{F}(z)$ defined on $G^{\ast}=\{z \in {\mathbb{C}}: \bar{z}
\in G\}$, we define a matrix-valued operator $\mathcal{F}(A^{\ast})$ as follows:
\begin{equation}
 \label{calc2}
<\mathcal{F}(A^{\ast})u^{\ast},v>=<u^{\ast}, \mathcal{F}^{T}(A)v>, ~~ \forall v \in
\mathcal{X}_{A, G},
\end{equation}
where $\mathcal{F}^T(A)$ is the matrix-valued operator with the symbol $\mathcal{F}(z)$ analytic in $G,$ and $u^{\ast} $ is an element of the space $\mathcal{X}^{'}_{A^{\ast}, G^{\ast}},$ dual to $\mathcal{X}_{A, G}.$ By construction, as dual to the space of direct product, the space  $\mathcal{X}^{'}_{A^{\ast}, G^{\ast}}$ represents the direct sum:
\[
\mathcal{X}^{'}_{A^{\ast}, G^{\ast}} = X^{'}_{A^{\ast}, G^{\ast}} \underset{\mbox{$m$ times}}{\underbrace{\oplus \dots \oplus}} X^{'}_{A^{\ast}, G^{\ast}},
\] 
with the corresponding topology. It follows from \eqref{mapping} that the mapping
\begin{equation}
\label{mapping2}
\mathcal{F}(A^{\ast}): \mathcal{X}^{'}_{A^{\ast}, G^{\ast}} \to \mathcal{X}^{'}_{A^{\ast}, G^{\ast}}
\end{equation}
is continuous. Indeed, assume that a sequence $u_n^{\ast} \in
\mathcal{X}^{'}_{A^{\ast}, G^{\ast}}$ converges
to $0$ in the topology of $\mathcal{X}^{'}_{A^{\ast}, G^{\ast}}.$ Then for arbitrary $v \in
\mathcal{X}_{A, G}$ we have
\[
<\mathcal{F}(A^{\ast}) u_n^{\ast}, v> = <u_n^{\ast}, \mathcal{F}^T(A)v>
=<u_n^{\ast},w>,
\]
where $w=\mathcal{F}^T(A)v \in \mathcal{X}_{A, G}$ due to \eqref{mapping}. 
Hence, $\mathcal{F}(A^{\ast})u_n^{\ast} \rightarrow 0,$
as $n \rightarrow \infty,$ in the topology of
$\mathcal{X}^{'}_{A^{\ast}, G^{\ast}},$ obtaining continuity of mapping \eqref{mapping2}.

\section{Main results}
\label{sec_main}

Below we derive representation formulas for solutions of fractional order systems of differential-operator equations.  We demonstrate the derivation in the case of Caputo fractional derivative. For the seek of clarity we start with the case $m=2$ and then the general case. The case of Riemann-Liouville fractional derivative can be treated similarly (see Section \ref{sec_RL}). 

\subsection{Fractional multi-order systems of differential-operator equations: $m=2$}

In this section we demonstrate the formal method of obtaining the representation formula for the solution of time-fractional arbitrary multi-order systems of differential-operator equations in the particular case of two equations. 
Namely, consider the system 
\begin{equation}
\label{200}
D_{\ast}^{\mathcal{B}}\mathcal{U}(t)=\mathbb{F}(A) \mathcal{U}(t) + \mathcal{H}(t),
\end{equation}
where $\mathcal{B}=(\beta_1, \beta_2), 0< \beta_1< \beta_2 \le 1,$ $\mathcal{H}(t)=(h_1(t), h_2(t))$ is a given vector-function, and  
\begin{equation}
\label{201}
\mathbb{F}(A) =
\begin{bmatrix}
f_{11}(A) & f_{1 2}(A) \\
f_{2 1}(A) & f_{2 2}(A)
\end{bmatrix},
\end{equation}
with the initial condition 
\begin{equation}
\label{ic200}
\mathcal{U}(0)=\varPhi =(\vf_1, \vf_2),
\end{equation} 
where $\varPhi \in  \mathcal{X}_{A,G}.$ We assume that $G$ does not contain roots of the equation 
\[
\Delta(z)={f_{1 1}(z)} f_{2 2}(z) - f_{2 1}(z) f_{1 2}(z)) = 0.
\] 
To find entries of the solution operator $\mathcal{S}(t,z)$ we consider the homogeneous counterpart of system \eqref{200} writing it in the explicit form
\begin{align*}
\begin{cases}
D_{\ast}^{\beta_1} u_1(t) &= f_{1 1}(A) u_1(t) + f_{1 2}(A) u_2(t),\\
D_{\ast}^{\beta_2} u_2(t) & = f_{2 1}(A) u_1(t) + f_{2 2}(A) u_2(t).
\end{cases}
\end{align*}
Applying the Laplace transform and replacing $A$ by the parameter $z$, we have
\begin{equation}
\label{lt200}
\begin{cases}
(Is^{\beta_1}-f_{1 1}(z) ) \mathcal{L}[u_1](s) - f_{1 2}(z) \mathcal{L}[u_2](s) &= Is^{\beta_1-1} \varphi_1 ,\\
- f_{2 1}(z) \mathcal{L}[u_1](s) + (Is^{\beta_2}-f_{2 2}(z) ) \mathcal{L}[u_2](s) &= Is^{\beta_2-1} \varphi_2.
\end{cases}
\end{equation}
The solution of system \eqref{lt200} is
\begin{equation}
\label{lts21}
\mathcal{L}[u_1](s) = \frac{1}{\Psi(s,z)} \Big(p_1(s,z) \vf_1 + q_{1}(s,z) \vf_2 \Big) \quad z \in G, s > r_{\ast}(z)
\end{equation}
\begin{equation}
\label{lts22}
\mathcal{L}[u_2](s) = \frac{1}{\Psi(s,z)} \Big( q_{2}(s,z) \vf_1+  p_2(s,z) \vf_2  \Big) , \quad z \in G, s > r_{\ast}(z).
\end{equation}
where 
\begin{align}
\label{lts23}
\Psi(s,z) &=s^{\beta_1+\beta_2}-s^{\beta_2}f_{1 1}(z)-s^{\beta_1}f_{2 2}(z)+\Delta(z),
\\
\label{lts24}
p_1(s,z) &=s^{\beta_1+\beta_2-1} -s^{\beta_1-1}f_{2 2}(z), \quad q_1(s,z) = s^{\beta_2-1}f_{1 2}(z),
\\
\label{lts25}
p_2(s,z) &=s^{\beta_1+\beta_2-1} -s^{\beta_2-1}f_{1 1}(z), \quad q_2(s,z) = s^{\beta_1-1}f_{2 1}(z),
\end{align}
and $r_{\ast}(z)$ is the real part of the roots of the equation $\Psi(s,z)=0.$ This solution is uniquely defined, since by assumption $G \cap Q = \emptyset,$ where $Q = \{z: \Psi(s,z)=0\}.$
We have
\begin{align}
\notag
\frac{1 }{\Psi(s,z)} &= \frac{1}{s^{\beta_2} \Big(s^{\beta_1} - f_{1 1}(z) -f_{2 2}(z)s^{\beta_1-\beta_2}+\Delta(z) s^{-\beta_2}\Big)}\\
&= \frac{1}{s^{\beta_2} \Big(s^{\beta_1}-f_{1 1}(z)\Big)\left( 1- \frac{f_{2 2}(z) s^{\beta_1-\beta_2}-\Delta(z) s^{-\beta_2}}{s^{\beta_1}-f_{1 1}(z)} \right) } \notag \\
&=\sum_{k=0}^{\infty}  \frac{ \left( s^{\beta_1-\beta_2}f_{2 2}(z)-s^{-\beta_2}\Delta(z) \right)^k}{s^{\beta_2} \left(s^{\beta_1}-f_{1 1}(z)\right)^{k+1}}
\notag \\
&= \sum_{k=0}^{\infty} \sum_{j=0}^k {k \choose j} f^j_{2 2}(z) \Big(-\Delta(z) \Big)^{k-j}  \frac{s^{ - k\beta_2+j\beta_1-\beta_2}}{ \Big(s^{\beta_1}-f_{1 1}(z)\Big)^{k+1}}.
\label{ser}
\end{align}
For $s>r_{\ast}(z)$ large enough the inequality 
\[
\left| \frac{s^{- k\beta_2+j\beta_1-\beta_2}}{ s^{\beta_1}-f_{1 1}(z)} \right| < 1
\]
is verified, and therefore, the series in \eqref{ser} is convergent. Now, for the solution of system \eqref{lt200}, we have
\begin{equation}
\label{ser01}
\mathcal{L}[u_1](s) = \sum_{k=0}^{\infty} \sum_{j=0}^k {k \choose j} f^j_{2 2}(z) \Big(-\Delta(z) \Big)^{k-j}  \frac{s^{ - k\beta_2+j\beta_1-\beta_2}}{ \Big(s^{\beta_1}-f_{1 1}(z)\Big)^{k+1}} \left[p_1(s,z) \vf_1 + q_{1}(s,z) \vf_2 \right],
\end{equation}
\begin{equation}
\label{ser02}
\mathcal{L}[u_2](s) = \sum_{k=0}^{\infty}  \sum_{j=0}^k {k \choose j} f^j_{2 2}(z) \Big(-\Delta(z) \Big)^{k-j}  \frac{s^{ - k\beta_2+j\beta_1-\beta_2}}{ \Big(s^{\beta_1}-f_{1 1}(z)\Big)^{k+1}}     \left[q_2(s,z) \vf_1 + p_{2}(s,z) \vf_2 \right].
\end{equation}
Consider the expressions
\[
P_1(s,z)= \frac{s^{ - k\beta_2+j\beta_1-\beta_2}p_1(s,z)}{ \Big(s^{\beta_1}-f_{1 1}(z)\Big)^{k+1}}, \quad Q_1(s,z)=  \frac{s^{ - k\beta_2+j\beta_1-\beta_2}q_1(s,z)}{ \Big(s^{\beta_1}-f_{1 1}(z)\Big)^{k+1}},
\]
and
\[
P_2(s,z)=  \frac{s^{ - k\beta_2+j\beta_1-\beta_2}p_2(s,z)}{ \Big(s^{\beta_1}-f_{1 1}(z)\Big)^{k+1}}, \quad Q_2(s,z)=  \frac{s^{ - k\beta_2+j\beta_1-\beta_2}q_2(s,z)}{ \Big(s^{\beta_1}-f_{1 1}(z)\Big)^{k+1}}.
\]
Further, let 
\[
\nu_{k j} =k \beta_2 - j \beta_1.
\] 
Since $\beta_2 > \beta_1$ and $k \ge j,$ we have $\nu_{k j} > 0,$ if $k \ge 1,$ and $\nu_{0 0}=0.$ Taking this into account, we obtain
\begin{align*}
P_1(s,z) 
&= \frac{1}{s^{\nu_{k j}}} \frac{s^{\beta_1 -1}}{ \Big(s^{\beta_1}-f_{1 1}(z)\Big)^{k+1}} -   \frac{f_{2 2}(z)}{s^{\nu_{k j} + \beta_2}} \frac{s^{\beta_1 -1}}{ \Big(s^{\beta_1}-f_{1 1}(z)\Big)^{k+1}}
\end{align*}
Now taking the inverse Laplace transform, due to formula \eqref{L2}, we obtain
\begin{align*}
\mathcal{L}^{-1}[P_1(\cdot,z)](t) &= J^{\nu_{k j}} \frac{t^{k \beta_1}}{k!} E^{(k)}_{\beta_1} (t^{\beta_1} f_{1 1}(z))
 -  f_{2 2} (z) J^{\nu_{k j}+\beta_2}  \frac{t^{k \beta_1}}{k!} E^{(k)}_{\beta_1} (t^{\beta_1} f_{1 1}(z)).
\end{align*}
Similarly,
\begin{align*}
\mathcal{L}^{-1}[Q_1(\cdot,z)](t) &= f_{1 2}(z) J^{\nu_{k j}+1} \frac{t^{k\beta_1+\beta_1-1} }{k!}E^{(k)}_{\beta_1, \beta_1} (t^{\beta_1} f_{1 1}(z)),
\end{align*}
\begin{align*}
\mathcal{L}^{-1}[P_2(\cdot,z)](t) &= J^{\nu_{k j}} \frac{t^{k \beta_1}}{k!}  E^{(k)}_{\beta_1} (t^{\beta_1} f_{1 1}(z))
 -   f_{1 1} (z) J^{\nu_{k j}+1} \frac{t^{k \beta_1+\beta_1-1}}{k!}   E^{(k)}_{\beta_1, \beta_1} (t^{\beta_1} f_{1 1}(z))  \Big].
\end{align*}
and
\begin{align*}
\mathcal{L}^{-1}[Q_2(\cdot,z)](t) &= f_{2,1}(z) J^{\nu_{k j}+\beta_2} \frac{t^{k \beta_1} }{k!}E^{(k)}_{\beta_1} (t^{\beta_1} f_{1 1}(z)),
\end{align*}
It follows from \eqref{ser01} and \eqref{ser02} that the entries $S_{jl}(t,z), j,l=1,2,$ of the matrix-symbol $\mathcal{S}(t,z)$ have representations:
\begin{align}
\notag
S_{1 1}(t,z) &= \sum_{k=0}^{\infty} \frac{1}{k!} \sum_{j=0}^k {k \choose j} f^j_{2 2}(z) \Big(-\Delta(z) \Big)^{k-j}     \Big[ J^{\nu_{k j}} \Big( t^{k \beta_1} E^{(k)}_{\beta_1} (t^{\beta_1} f_{1 1}(z)) \Big)
\\
& -  f_{2 2} (z) J^{\nu_{k j}+\beta_2} \Big(  t^{k \beta_1} E^{(k)}_{\beta_1} (t^{\beta_1} f_{1 1}(z)) \Big) \Big] \notag
\\
&= \Big(I - f_{2 2}(z)J^{\beta_2}\Big)\sum_{k=0}^{\infty} { \Big( f_{2 2}(z) J^{\beta_2-\beta_1} - \Delta(z)J^{\beta_2} \Big)^k} \Big[\frac{ t^{k \beta_1}}{k!}  E^{(k)}_{\beta_1}(t^{\beta_1} f_{1 1}(z)) \Big],
\label{sol200}
\end{align}

\begin{align}
S_{1 2}(t,z) &=  f_{1 2}(z) \sum_{k=0}^{\infty}  \frac{1}{k!}\sum_{j=0}^k {k \choose j} f^j_{2 2}(z) \Big(-\Delta(z) \Big)^{k-j}   J^{\nu_{k j}+1} \Big( t^{k\beta_1+\beta_1-1}  E^{(k)}_{\beta_1, \beta_1} (t^{\beta_1} f_{1 1}(z)) \Big) \notag
\\
&= f_{1 2}(z) \sum_{k=0}^{\infty}  { \Big( f_{2 2}(z) J^{\beta_2-\beta_1} - \Delta(z)J^{\beta_2} \Big)^k}  J \Big( \frac{t^{k\beta_1+\beta_1-1}}{k!}   E^{(k)}_{\beta_1, \beta_1} (t^{\beta_1} f_{1 1}(z)) \Big),
\label{sol201}
\end{align}
\begin{align}
S_{2 1}(t,z) &=  f_{2 1}(z)\sum_{k=0}^{\infty} \frac{1}{k!}\sum_{j=0}^k {k \choose j} f^j_{2 2}(z) \Big(-\Delta(z) \Big)^{k-j}
  J^{\nu_{k j}+\beta_2} \Big( t^{k \beta_1}  E^{(k)}_{\beta_1} (t^{\beta_1} f_{1 1}(z)) \Big) \notag
  \\
  &= f_{2 1}(z)\sum_{k=0}^{\infty} { \Big( f_{2 2}(z) J^{\beta_2-\beta_1} - \Delta(z)J^{\beta_2} \Big)^k}
  J^{\beta_2} \Big(\frac{ t^{k \beta_1}}{k!}   E^{(k)}_{\beta_1} (t^{\beta_1} f_{1 1}(z)) \Big),
 \label{sol202}
\\
S_{2 2}(t,z) &= \sum_{k=0}^{\infty} \frac{1}{k!} \sum_{j=0}^k {k \choose j} f^j_{2 2}(z) \Big(-\Delta(z) \Big)^{k-j}
 \Big[  J^{\nu_{k j}} \Big( t^{k \beta_1}  E^{(k)}_{\beta_1} (t^{\beta_1} f_{1 1}(z)) \Big)
 \notag
\\
& -  f_{1 1} (z) J^{\nu_{k j}+1} \Big( t^{k \beta_1+\beta_1-1}   E^{(k)}_{\beta_1, \beta_1} (t^{\beta_1} f_{1 1}(z)) \Big) \Big], \notag
\\
&= \sum_{k=0}^{\infty} { \Big( f_{2 2}(z) J^{\beta_2-\beta_1} - \Delta(z)J^{\beta_2} \Big)^k}
 \Big[ \frac{ t^{k \beta_1}}{k!}   E^{(k)}_{\beta_1} (t^{\beta_1} f_{1 1}(z))
 \notag
\\
& -  f_{1 1} (z) J \Big(\frac{ t^{k \beta_1+\beta_1-1}}{k!}    E^{(k)}_{\beta_1, \beta_1} (t^{\beta_1} f_{1 1}(z)) \Big) \Big] .
\label{sol203}
\end{align}

Thus, we proved the following theorem.
\begin{thm} \label{thm_m=20}
The formal solution to system \eqref{200}, \eqref{ic200} has the representation 
\[
\mathcal{U}(t) = \mathcal{S}(t,A) \Phi + \int_0^{t} \mathcal{S}(t-\tau,A) D_{+}^{1-\mathcal{B}} \mathcal{H}(\tau) d\tau
\]
where $\mathcal{S}(t,A)$ is the matrix-valued solution operator with the matrix-symbol $\mathcal{S}(t,z),$ entries of which are defined in \eqref{sol200}-\eqref{sol203}.
\end{thm}

\begin{thm}
\label{thm_m=21}
Let $f_{1 2}(z) = 0.$ Then the solution operator $\mathcal{S}(t,A)$ has the matrix-symbol with entries
\begin{align}
\label{m=11}
S_{1 1}(t,z) & =E_{\beta_1}(t^{\beta_1}f_{1 1}(z)), \quad S_{1 2}(t,z)=0,  
\\
\label{m=21}
 S_{2 1}(t,z) &= f_{2 1}(z) \Big( E_{\beta_1}(t^{\beta_1}f_{1 1}(z)) \Big) \ast \Big( t^{\beta_2-1} E_{\beta_2, \beta_2}(t^{\beta_2}f_{2 2}(z)) \Big),
\\
\label{m=22}
 S_{2 2} (t, z) &= E_{\beta_2}(t^{\beta_2}f_{2 2}(z)),
 \end{align}
 where $"\ast"$ is the convolution operation.
\end{thm}
{\it Proof} The fact that $S_{1 2}(t,z)=0$ obviously follows from \eqref{sol201}. Now we show the equality for $S_{1 1}(t,z).$ 
First we notice  that $f_{1 2}(z)=0$ implies $\Delta(z)=f_{1 1}(z) f_{2 2}(z).$ Taking this fact into account and utilizing the semigroup property \cite{SKM93,Umarov_book2015} of the fractional integration operator $J^{\beta}$, we can express $S_{1 1}(t,z)$ in the form
\begin{align}
S_{1 1}(t,z) &= \Big(I- f_{2 2}(z)J^{\beta_2}\Big)\sum_{k=0}^{\infty}  \Big( f_{2 2}(z) J^{\beta_2-\beta_1} - f_{1 1}(z) f_{2 2}(z)J^{\beta_2} \Big)^k
 \Big[ \frac{ t^{k \beta_1}}{k!} E^{(k)}_{\beta_1}(t^{\beta_1} f_{1 1}(z)) \Big]  \notag 
\\
&= \Big(I- f_{2 2}(z)J^{\beta_2}\Big) \sum_{k=0}^{\infty}  \Big( f_{2 2}(z) J^{\beta_2-\beta_1} \Big)^k  
 \Big( I-f_{1 1}J^{\beta_1} \Big)^k 
\Big[ \frac{t^{k \beta_1}}{k!} E^{(k)}_{\beta_1}(t^{\beta_1} f_{1 1}(z)) \Big],
\label{lower}
\end{align}
where $I$ is the identity operator. Now,  due to Lemma \ref{lemma_1} with $\nu=1$ and $\mu=f_{1 1}(z)$, one has
\begin{equation}
\label{rel_new0}
\Big( I-f_{1 1}(z)J^{\beta_1} \Big)^k\Big[ \frac{t^{k \beta_1}}{k!} E^{(k)}_{\beta_1}(t^{\beta_1} f_{1 1}(z)) \Big] = J^{k \beta_1} E_{\beta_1}(t^{\beta_1}f_{1 1}(z)),
\end{equation}
valid for all $k = 0,1,\dots.$ Thus \eqref{lower} reduces to
\begin{align}
S_{1 1}(t,z) 
&=\Big(I- f_{2 2}(z)J^{\beta_2}\Big)  \sum_{k=0}^{\infty}  \Big( f_{2 2}(z) J^{\beta_2-\beta_1} \Big)^k   J^{k \beta_1} E_{\beta_1}(t^{\beta_1}f_{1 1}(z)) \notag 
\\
&=\Big(I- f_{2 2}(z)J^{\beta_2}\Big) \sum_{k=0}^{\infty}  \Big( f_{2 2}(z) J^{\beta_2-\beta_1}  J^{\beta_1} \Big)^k    
E_{\beta_1}(t^{\beta_1}f_{1 1}(z)),
\label{lower1}
\end{align}
Since
\[
\sum_{k=0}^{\infty}  \Big( f_{2 2}(z) J^{\beta_2} \Big)^k = \Big( I - f_{2 2}(z) J^{\beta_2} \Big)^{-1},
\]
it follows from \eqref{lower1} the desired equality for $S_{1 1}(t,z)$ in \eqref{m=11}.

Now we show the validity of \eqref{m=21}. Taking into account $\Delta(z)=f_{1 1}(z) f_{2 2}(z)$, we rewrite $S_{2 1}(t,z)$ in \eqref{sol202} in the form
\begin{align*}
S_{2 1}(t,z)
  &= f_{2 1}(z) \sum_{k=0}^{\infty}  \Big( f_{2 2}(z) J^{\beta_2-\beta_1} - f_{1 1}(z) f_{2 2}(z) J^{\beta_2} \Big)^k
  J^{\beta_2} \Big( \frac{ t^{k \beta_1}}{k!}  E^{(k)}_{\beta_1} (t^{\beta_1} f_{1 1}(z)) \Big)
  \\
  &= f_{2 1}(z)\sum_{k=0}^{\infty}   \Big( f_{2 2}(z) J^{\beta_2-\beta_1} \Big)^k  
J^{\beta_2} 
\Big( I-f_{1 1}J^{\beta_1} \Big)^k  \Big(\frac{ t^{k \beta_1}}{k!}  E^{(k)}_{\beta_1} (t^{\beta_1} f_{1 1}(z)) \Big) ,
 \label{m=21_1}
\end{align*}
Using the relation \eqref{rel_new0}, we have
\begin{align*}
S_{2 1}(t,z)
= & f_{2 1}(z) \Big( I - f_{2 2}(z) J^{\beta_2} \Big)^{-1} J^{\beta_2} E_{\beta_1}(t^{\beta_1}f_{1 1}(z)).
  \end{align*}
Further, due to Lemma \ref{lemma_2} (see Eq. \eqref{rel_new3}) with $\mu_1=f_{1 1}(z)$ and $\mu_2=f_{2 2}(z),$ we obtain
relation \eqref{m=21}.

Similarly, $S_{2 2}(t,z)$ in \eqref{sol203} can be written as 
\begin{align*}
S_{2 2}(t,z) 
&= \sum_{k=0}^{\infty} \Big( f_{2 2}(z) J^{\beta_2-\beta_1} - \Delta(z)J^{\beta_2} \Big)^k
 \Big[  \frac{ t^{k \beta_1}}{k!}  E^{(k)}_{\beta_1} (t^{\beta_1} f_{1 1}(z))
 \notag
\\
& -  f_{1 1} (z) J \Big( \frac{t^{k \beta_1}}{k!}   E^{(k)}_{\beta_1, \beta_1} (t^{\beta_1} f_{1 1}(z)) \Big) \Big] 
\\
&= \sum_{k=0}^{\infty}   \Big( f_{2 2}(z) J^{\beta_2-\beta_1} \Big)^k  \Big( I-f_{1 1}J^{\beta_1} \Big)^k  \Big( \frac{ t^{k \beta_1}}{k!}  E^{(k)}_{\beta_1} (t^{\beta_1} f_{1 1}(z)) 
\notag
\\
& - f_{1 1}(z) J \frac{ t^{k \beta_1+\beta_1-1}}{k!}  E^{(k)}_{\beta_1, \beta_1} (t^{\beta_1} f_{1 1}(z)) \Big).
\end{align*}
Finally, using Lemma \ref{lemma_1}, we obtain
\begin{align*}
S_{2 2}(t,z) &= \Big( I - f_{2 2}(z) J^{\beta_2} \Big)^{-1} \Big( E_{\beta_1}(t^{\beta_1}f_{1 1}(z)) - f_{1 1}(z) J E_{\beta_1, \beta_1}(t^{\beta_1}f_{1 1}(z)) \Big)
\\
&= E_{\beta_2} (t^{\beta_2}f_{2 2}(z)).
\end{align*}
In the last step we used relation \eqref{rel_new3} with $\mu_1=f_{1 1}(z)$ and $\mu_2=f_{2 2}(z).$
\eproof
\begin{rem}
Theorem \ref{thm_m=21} states that the representation formula presented in Theorem \ref{thm_m=20} coincides with the representation formula obtained in \cite{Umarov2024} for the solution of fractional order systems with a lower triangular matrix-valued operator.
\end{rem}

\subsection{Fractional multi-order systems of differential-operator equations: $m \ge 2$}
\label{sec_3.2}
The method demonstrated in the previous section for $m=2$ works, in fact, for arbitrary number of equations. To derive the solution operator consider the homogeneous system
\begin{equation}
\label{GEN-1}
D_{\ast}^{\mathcal{B}}\mathcal{U}(t)=\mathcal{F}(A) \mathcal{U}(t), \quad t>0,
\end{equation}
with the initial condition 
\begin{equation}
\label{GEN-2}
\mathcal{U}(0)= \varPhi \in  \mathcal{X}_{A,G}.
\end{equation}
Here $\mathcal{B}=(\beta_1,\dots,\beta_m), 0< \beta_j \le 1, \ j=1,\dots, m,$ is arbitrary vector-order. We can assume without loss of generality that $\beta_1 = \min \{\beta_1,\dots,\beta_m\}.$
Applying the Laplace transform and replacing $A$ by a parameter $z \in G$, we obtain
\begin{equation}
\label{am03}
(Is^{\mathcal{B}} - \mathcal{F}(z)) \mathcal{L} [ \mathcal{U}] (s) = Is^{\mathcal{B}-{\bf 1}} \varPhi,
\end{equation}
where $Is^{\mathcal{B}}=\text{diag}(s^{\beta_1}, \dots, s^{\beta_m}).$ This is a system of linear algebraic equations dependent on parameters $s \in \mathbb{C}$ and $z \in G.$  The determinant of the matrix on the left has the structure
\[
\Psi(s,z)=\det(Is^{\mathcal{B}}-\mathcal{F}(z))= s^{|\beta|} + \mathcal{G}(s,z) + (-1)^m \Delta(z),
\]
where $|\beta|=\beta_1+\dots + \beta_m,$ $\Delta(z) = \det{\mathcal{F}(z)}$ and $\mathcal{G}(s,z)$ has the form
\begin{equation}
\label{g00}
\mathcal{G}(s,z) = g_{1}(z)s^{\beta_1} +\dots + g_m(z) s^{\beta_m} + g_{1 2} s^{\beta_1+\beta_2} + g_{1 3} s^{\beta_1+\beta_3} + \dots + g(z) s^{\beta_2+\dots+\beta_m}.
\end{equation}
In other words $\mathcal{G}(s,z)$ is the sum of functions of the form
\[
\sum_{\alpha} g_{\alpha}(z) s^{\beta_{\alpha_p}+\dots + \beta_{\alpha_q}},
\]
where $\alpha=(\alpha_p,\dots,\alpha_q), \ 1\le p < q \le m,$ is a multi-index taking values in the subsets of the set  $\{1,\dots,m\}$ except  $(1,\dots,m),$ the function $g_{\alpha}(z)$ is the sum and multiplication combinations of entrees $f_{kj}(z)$ of the matrix-symbol $\mathcal{F}(z).$ 
Let $r_{\ast}$ is the real part of the largest root of the equation $\Psi(s,z)=0.$ Then, for $s>r_{\ast}$ system \eqref{am03} has a unique solution
\begin{equation}
\label{LSOLUTION}
 \mathcal{L} [ \mathcal{U}] (s) = (Is^{\mathcal{B}} - \mathcal{F}(z))^{-1} Is^{\mathcal{B}-{\bf 1}} \varPhi,
\end{equation}
It follows that the solution operator has the matrix-valued symbol
\begin{equation}
\label{SOLUTION_SYMBOL}
\mathcal{S}(t,z) = \mathcal{L}^{-1}_{s \to t} \Big[(Is^{\mathcal{B}} - \mathcal{F}(z))^{-1} Is^{\mathcal{B}-{\bf 1}}\Big].
\end{equation}

The components of the solution have the structure (as an implication of the well-known Cramer's rule)
\[
\mathcal{L}[{u}_j](s) =\frac{\mathcal{P}_{ j}(s,z)}{\Psi(s,z)}, \quad j=1,\dots,m,
\]
where $u_j=u_j(t)$ is the $j$-th component of $\mathcal{U}(t)$ and $\mathcal{P}_j(s,z)$ is the determinant of the matrix obtained by replacing $j$-th column of the matrix $Is^{\mathcal{B}}-\mathbb{F}(z)$ with the vector $Is^{\mathcal{B}-{\bf 1}} \varPhi.$ 
The latter can be rewritten in the form
\begin{equation}
\label{solm00}
\mathcal{L}[{u}_j](s) =\sum\limits_{l=1}^m \frac{{p}_{j l}(s,z) }{\Psi(s,z)} \vf_l,  \quad j=1,\dots,m,
\end{equation}
where $\vf_l, l=1,\dots,m,$ are components of $\varPhi,$ and functions ${p}_{j l}(s,z)$ have form \eqref{g00}.
Let $\beta_{\ast}=\beta_2+\dots + \beta_m.$ We have
\[
\frac{1 }{\Psi(s,z)} = \frac{1}{s^{\beta_{\ast}}\Big(s^{\beta_1} + g(z) - g(s,z) s^{-\beta_{\ast}} +(-1)^m\Delta(z) s^{-\beta_{\ast}} \Big)}
\]
where $g(z)$ is the coefficient in the term $g(z)s^{\beta_2+\dots+\beta_m}$ in \eqref{g00} and 
\[
g(s,z)=-G(s,z)+g(z)s^{\beta_2+\dots+\beta_m}.
\]
Further, similar to the case $m=2,$ we represent $1/\Psi(s,z)$ in the infinite functional series form
\begin{align}
\label{}
\frac{1}{\Psi(s,z)}&=
\frac{1}{s^{\beta_{\ast}} \Big(s^{\beta_1}+g(z)\Big)\left( 1- \frac{g(s,z)s^{-\beta_{\ast}}+(-1)^m \Delta(z)s^{-\beta_{\ast}}}{s^{\beta_1}+g(z)} \right) } \notag \\
&=\sum_{k=0}^{\infty}  \frac{  \Big( g(s,z) + (-1)^m \Delta(z) \Big)^k s^{-k \beta_{\ast} -\beta_{\ast}} }{ \Big(s^{\beta_1}+g(z)\Big)^{k+1}}.
\label{serm}
\end{align}
For $s>r_{\ast}(z)$ large enough the inequality 
\[
\left| \frac{g(s,z)s^{-\beta_{\ast}}+(-1)^m \Delta(z)s^{-\beta_{\ast}}}{s^{\beta_1}+g(z)} \right| < 1
\]
is valid, and therefore, the series in \eqref{serm} is convergent.

Further, it follows from \eqref{solm00} that 
\[
u_j(t)= \sum_{l=0}^{m} \mathcal{L}^{-1} \left[ \frac{{p}_{j l}(s,A) }{\Psi(s,z)} \right] \vf_j, \quad j=1,\dots,m.
\]
Hence, the matrix-symbol of the solution-operator $\mathcal{S}(t,A)$ has entries
\[
s_{j l} (t, z) = \mathcal{L}^{-1}_{s \to t} \left[ \frac{{p}_{j l}(s, z) }{\Psi(s, z)} \right](t), \quad j, l = 1, \dots, m.
\]
Now using \eqref{serm} we have
\[
s_{j l}(s,z)= \mathcal{L} [s_{j l}(t,z)](s) =\sum_{k=0}^{\infty}  \frac{p_{j l}(s,z)  \Big( g(s,z) + (-1)^m \Delta(z) \Big)^k s^{-k \beta_{\ast} -\beta_{\ast}} }{ \Big(s^{\beta_1}+g(z)\Big)^{k+1}}.
\]
Taking into account the fact that ${g}(s,z)$ and functions $p_{j l}(s,z)$ can be represented in form \eqref{g00}, we can write the expressions $p_{j l}(s,z) ({g}(s,z)+(-1)^m \Delta(z))^k, \ j,l=1,\dots,m,$ as
\begin{equation}
\label{q}
\sum\limits_{\alpha} Q_{\alpha, j, l, k}(z) s^{\gamma_{\alpha, j,l,k}},
\end{equation}
where $Q_{\alpha, j, l, k}(z)$ are sum and product combinations of $f_{j,k}(z)$ and exponents $\gamma_{\beta, j, l, k}$ depend on the sum combinations of $\beta_1, \dots, \beta_m$ and their multiples.
Therefore,
\begin{align*}
s_{j l}(s,z) =\sum_{k=0}^{\infty} \sum\limits_{\alpha} Q_{\alpha, j, l, k}(z)\frac{s^{\gamma_{\alpha, j,l,k} - (k+1) \beta_{\ast}} }{\Big(s^{\beta_1}+g(z)\Big)^{k+1}}.
\end{align*}
Further, let $\nu_{\alpha, j, l, k} = (k+1)\beta_{\ast} + \beta_1-1 - \gamma_{\alpha,j,l,k}.$ By construction $\nu_{\alpha, j, l, k} \ge 0$ for all indices $\alpha, j, l, k.$ Then,
\begin{align}
s_{j l}(t,z) &= \sum_{k=0}^{\infty} \sum\limits_{\alpha} Q_{\alpha, j, l, k}(z) \mathcal{L}_{s \to t}^{-1} \left[ \frac{s^{\beta_1 - 1}}{s^{\nu_{\alpha, j, l, k}}\Big(s^{\beta_1}+g(z)\Big)^{k+1}}\right] \notag
\\
&= \sum_{k=0}^{\infty} \sum\limits_{\alpha} Q_{\alpha, j, l, k}(z) J^{\nu_{\alpha, j, l, k}} \left[ \frac{t^{k\beta_1}}{k!}  E^{(k)}_{\beta_1} \Big(-t^{\beta_1}g(z) \Big) \right], \quad j, l =1,\dots,m.
\label{solm5}
\end{align}
Hence, the solution to Cauchy problem \eqref{GEN-1}, \eqref{GEN-2} has the form
\begin{equation}
\label{GEN-3}
\mathcal{U}(t) = \mathcal{S}(t,A) \varPhi,
\end{equation}
where $\mathcal{S}(t,A)$ is the matrix-valued solution operator with the matrix-symbol $\mathcal{S}(t,z)$ defined by \eqref{solm5}.

\begin{thm} 
\label{thm_001}
Let $\mathcal{B}= ( \beta_1,\dots,\beta_m ),$ where $\beta_j \in (0,1), j=1,\dots, m,$ arbitrary numbers and $\beta_1=\min\{\beta_1,\dots,\beta_m\}$. Let 
$A$ be a closed operator defined on a Banach space $X,$ the set $G$ satisfies the condition $G \cap Q_0 = \emptyset,$ where $Q_0=\{z: \Psi(s,z) =0\}$, $ \varPhi \in \mathcal{X}_{A,G},$    $H(t) \in AC[\re_+; \mathcal{X}_{A,G}],$ and $D_+^{1-\mathcal{B}}H(\tau,x) \in C[\re_+;  \mathcal{X}_{A,G} ].$   

Then for any $T>0$ Cauchy problem \eqref{01}, \eqref{02} 
has a unique solution $U(t)\in C^{\infty}[(0,T]; \mathcal{X}_{A,G}] \cap C[[0,T]; \mathcal{X}_{A,G}],$ having the representation
\begin{equation}
\label{solution_1}
U(t) = S(t,A)\varPhi + \int\limits_0^t S(t-\tau, A) D_+^{1-\mathcal{B}}H(\tau) d\tau, \quad t>0, 
\end{equation}
where $S(t,A)$ is the solution matrix-operator with the matrix-symbol $\mathcal{S}(t, z)$ defined in \eqref{solm5}.
\end{thm}

\textit{Proof.} The representation \eqref{solution_1} follows  from \eqref{GEN-3} and fractional Duhamel's principle \cite{Umarov_book2015,Umarov_handbook}. Let
\begin{align}
\label{Cauchy_01000_h}
V(t) &= S(t,A)\varPhi, \\ 
W(t) &= \int\limits_0^t S(t-\tau, A) D_+^{1-\mathcal{B}}H(\tau) d\tau, \quad t \ge 0.
\label{Cauchy_02000_h}
\end{align}
It follows from \eqref{mapping} that $V(t) \in { \mathcal{X}_{A, G}}$ for every fixed $t \ge 0,$ continuous on $[0,T],$ and infinitely differentiable on $(0,T)$ in the topology of $\mathcal{X}_{A, G}$ due to the construction of the solution operator $S(t,A).$ 
Similarly, by virtue of the continuity in $\mathcal{X}_{A,G}$ of operators with symbols analytic in $G$  (see \eqref{mapping}),  for each fixed $t$ we have $W(t) \in  \mathcal{X}_{A, G}$ 
for each fixed $t \in [0,T].$ The continuity of $W(t)$ on $[0,T]$ in the variable $t$ and its infinite differentiability on $(0,T)$ follows from the construction of the solution operator $S(t,D)$ in the standard way.
\eproof

\begin{thm}
\label{thm_multi_1_dual}
Let $X$ be a reflexive Banach space with the conjugate $X^{\ast},$ $A$ be a closed operator with a domain $\mathcal{D} \subset X,$ and $\mathcal{F}(A)$ be a matrix operator with the symbol $\mathcal{F}(z)$ continuous on $G.$ and satisfying condition $G \cup Q_0 = \emptyset,$ where $Q_0$ is defined in \eqref{Q0}. Assume that   $\varPsi \in \mathcal{X}^{\prime}_{A^{\ast}, G^{\ast}},$    $H(t) \in AC[\re_+; \mathcal{X}^{\prime}_{A^{\ast}, G^{\ast}}],$ and $D_+^{1-\mathcal{B}}H(t) \in C[\re_+;  \mathcal{X}^{\prime}_{A^{\ast}, G^{\ast}} ].$   Then for any $T>0$ Cauchy problem
\begin{align} \label{Cauchy_01_dual}
D_{\ast}^{\mathcal{B}} {V}(t) &= \mathcal{F}(A^{\ast}) {V} (t) + H(t), \quad t>0, \\
V(0) & =\varPsi, 
\label{Cauchy_02_dual}
\end{align}
has a unique solution $V(t) \in C^{\infty}[(0,T]; \mathcal{X}^{\prime}_{A^{\ast}, G^{\ast}}] \cap C[[0,T]; \mathcal{X}^{\prime}_{A^{\ast}, G^{\ast}}],$ having the representation
\begin{equation}
\label{solution_1_dual}
V(t) = S(t,A^{\ast})\varPsi +\int\limits_0^t S(t-\tau, A^{\ast}) D_+^{1-\mathcal{B}}H(\tau) d\tau, \quad t>0, 
\end{equation}
where $S(t,A^{\ast})$ is the operator with the matrix-symbol $\mathcal{S}(t, z)$ defined in \eqref{sol_0001}. 
\end{thm}

\textit{Proof.} We note that elements $ D_{\ast}^{\mathcal{B}} {V}(t) $ and $ \mathcal{F}(A^{\ast}) {V} (t) $ belong to the space $\mathcal{X}^{\prime}_{A^{\ast}, G^{\ast}}$ if $V(t) \in \mathcal{X}^{\prime}_{A^{\ast}, G^{\ast}}$ for each fixed $t \ge 0.$ This fact follows from the definition of the fractional derivative $D_{\ast}^{\mathcal{B}}$ and Theorem \ref{thm_001}.

We show that $V(t)$ defined in \eqref{solution_1_dual} 
 satisfies the following conditions:
\begin{align} \label{Cauchy_01_dual_1}
\langle D_{\ast}^{\mathcal{B}} {V}(t) , \varPhi \rangle &=  \langle {V} (t) , \mathcal{F}^{T}(A) \varPhi \rangle + \langle H(t), \varPhi \rangle, \quad t>0, 
\\
\langle V(0), \varPhi \rangle & = \langle \varPsi(x), \varPhi \rangle, 
\label{Cauchy_02_dual_1}
\end{align}
where $\mathcal{F}^{T}(A)$ is the conjugate transpose of $\mathcal{F}(A^{\ast}),$ for an arbitrary element $\varPhi$ in the space $\mathcal{X}_{A, G}.$ 
Indeed, to show this fact let us first assume that $H(t) = 0$\footnote{as an element of $\mathcal{X}^{\prime}_{A^{\ast}, G^{\ast}}$} for all $t \ge 0.$ Then \eqref{Cauchy_01_dual_1} takes the form
\begin{equation}
\left\langle \Big[ D_{\ast}^{\mathcal{B}} {S}(t,A^{\ast}) - \mathcal{F}(A^{\ast}) \Big] V(t), \varPhi \right\rangle =  
\left\langle V(t), \Big[ D_{\ast}^{\mathcal{B}} {S}(t,A) - \mathcal{F}(A) \Big]^T \varPhi \right\rangle = 0, 
\end{equation}
for each fixed $t>0.$ The operator $S(t,A)$ is constructed so that $D_{\ast}^{\mathcal{B}} {S}(t,A) - \mathcal{F}(A)=0,$ which implies 
$[D_{\ast}^{\mathcal{B}} {S}(t,A) - \mathcal{F}(A)]^T=0,$ as well. Indeed, if $V(t)$ is a solution to equation \eqref{GEN-1}, then it follows from representation \eqref{GEN-3} that $D_{\ast}^{\mathcal{B}} V(t) = D_{\ast}^{\mathcal{B}} {S}(t,A) \varPhi = \mathcal{F}(A) \varPhi $ for any fixed $\varPhi \in \mathcal{X}_{A, G}.$ This implies the equality $D_{\ast}^{\mathcal{B}} {S}(t,A) = \mathcal{F}(A).$ Thus, equation \eqref{Cauchy_01_dual_1} is valid for all $\varPhi \in \mathcal{X}_{A, G}$. 

Further, it follows from the construction of the operator $\mathcal{S}(t,A)$ that the symbol $\mathcal{S}(t,z)$ at $t=0$ reduces to the identity matrix.  Therefore, the operator corresponding to the matrix-symbol $\mathcal{S}(0,z)$ is the identity operator. Hence, $V(0) = S(0,A^{\ast}) \varPhi= \varPhi.$ Thus, equality \eqref{Cauchy_02_dual_1} is also verified.

In the general case, for non-zero $H(t),$ the representation \eqref{solution_1_dual} is an implication of the fractional Duhamel principle \cite{Umarov_book2015,Umarov_handbook}.
\eproof

\begin{rem}
Obtaining closed form representations for $Q_{\alpha,j,l,k}$ and $\nu_{\alpha,j,l,k}$ in the general case is possible, but it is very cumbersome. In the case of $m=2$ see \eqref{sol200}--\eqref{sol203}. 
\end{rem}

\subsection{Fractional multi-order systems of differential-operator equations with triangular matrix-valued operators}

If  the matrix-symbol $\mathcal{F}(z)$ in system \eqref{01} is lower or upper triangular matrix then representation \eqref{solm5} significantly simplifies. See Theorem \ref{thm_m=21} in the case of lower triangular matrix for $m=2.$ In this section we derive representation formulas for the solution for arbitrary $m \ge 2.$ For fractional systems of pseudo-differential equations with lower or upper triangular matrix-symbols the representation formulas are presented in paper \cite{Umarov2024}.

Assume that $\mathcal{B}$ is an arbitrary multi-order with components $\beta_k \in (0,1), \ k=1,\dots,m$ and $\mathcal{F}(z)$ is a lower triangular matrix-symbol. Then system \eqref{am03} takes the form 
\[
\begin{bmatrix}
s^{\beta_1}-f_{1 1}(z)  & 0 & \dots & 0 \\
-f_{2 1} (z) & s^{\beta_2}-f_{2 2}(z) & \dots & 0 \\
\dots & \dots & \dots & \dots \\
-f_{m 1}(z) & -f_{m 2}(z) & \dots & s^{\beta_m}-f_{m m}(z)
\end{bmatrix} 
\begin{bmatrix}
L[u_1](s) \\
L[u_2](s) \\
\dots \\
L[u_m](s)
\end{bmatrix}
=
\begin{bmatrix}
s^{\beta_1-1} \vf_1 \\
s^{\beta_2-1} \vf_2 \\
\dots\\
s^{\beta_m-1} \vf_m
\end{bmatrix}.
\]
The latter implies the following recurrent relations
\begin{align}
L[u_1](s) &= \frac{s^{\beta_1-1}}{s^{\beta_1}-f_{1 1}(z)} \vf_1,
\label{case4_1} \\
L[u_k](s) & = \frac{s^{\beta_k-1}}{s^{\beta_k}-f_{1 1}(z)} \vf_k + \sum_{j=1}^{k-1} \frac{f_{k j}(z)}{s^{\beta_k}-f_{k k}(z)} L[u_j](s), \quad k=2,\dots,m.
\label{case4_2}
\end{align}

In order to represent the solution \eqref{case4_2} through $\vf_1, \dots, \vf_m,$ we introduce the following notations. Let $\mathcal{T}=\{k,k-1,\dots,k-j \},$ where $k$ and $j$  are integers satisfying conditions $1 \le k \le m$ and $1\le j < k,$ respectively. By $\mathcal{T}_l, \ 1 \le l \le j -1,$ we denote the set of subsets of $\mathcal{T}$ such that from $\mathcal{T}$ exactly $l$ numbers except $k$ and $k-j$ are removed. Then for $L[u_k](s), \ k=2,\dots, m,$ in \eqref{case4_2}, we have

\begin{equation}
\label{case4_3}
L[u_k](s) = \frac{s^{\beta_k-1}}{s^{\beta_k}-f_{1 1}(z)} \vf_k + \sum_{j=1}^{k-1}  \sum_{l=1}^{j-1} \frac{P_{k j l}(z) s^{\beta_{k-j}-1}}{\prod\limits_{\tau \in \mathcal{T}_l} \Big(s^{\beta_{\tau}}-f_{\tau \tau}(z)\Big)} \vf_{k-j}, \quad k=2,\dots,m,
\end{equation}
where $P_{k j l}, \ k=2,\dots, m, \ j= 1,\dots,k-1,$ are multiplication and sum combinations of functions $f_{\tau \mu}, \ \tau=2,\dots, k, \ \mu=1,\dots,\nu-1.$
Now, making use of formula \eqref{L1} and the convolution formula for the Laplace transform, it follows from \eqref{case4_1} and \eqref{case4_3} that
\begin{align}
&u_1(t) = E_{\beta_1}(f_{1 1}(A) t^{\beta_1})\vf_1,\\
&u_k(t) = E_{\beta_k}(f_{k k}(A)t^{\beta_k}) \vf_k  \notag \\
&+ \sum_{j=1}^{k-1} \sum_{l=1}^{j-1} \left[ P_{k j l}(A) E_{\beta_{k-j}} \Big(f_{k-j k-j}(A)t^{\beta_{k-j}}\Big) \ast \Big( \ast \prod_{\tiny \begin{matrix} \tau \in \mathcal{T}_l \\ \tau \neq j \end{matrix}}  t^{\beta_{\tau}-1}E_{\beta_{\tau}, \beta_{\tau}}(f_{\tau \tau}(A)t^{\beta_{\tau}})  \Big)\right] \vf_{k-j},
\\
& \hspace{6cm} k=2,\dots,m, \notag
\end{align}
where $J^{\gamma}$ is the fractional integration operator (see \eqref{FrInt_1}) of order $\gamma > 0,$ "$\ast$" is the convolution operation, and "$\ast \prod $" is the convolution product. Thus, the solution of homogeneous Cauchy problem \eqref{GEN-1}, \eqref{GEN-2} has the representation
\begin{equation}
\label{case4_S1}
\mathcal{U}(t) = \mathcal{S}(t,A) \varPhi,
\end{equation}
where $\mathcal{S}(t,A)$ is the solution matrix-operator with the matrix-symbol $\mathcal{S}(t,z)$ with entries
\begin{equation}
s_{k \, j}(t,z) = 
\begin{cases}
0, & \text{if }  j >k, \\
E_{\beta_k}(f_{k k}(z)t^{\beta_k}), & \text{if } j=k, \\
\sum\limits_{l=1}^{j-1} P_{k j l}(z) E_{\beta_j}(f_{j j}(z)t^{\beta_j}) \ast \Big( \ast \prod\limits_{\tiny \begin{matrix} \tau \in \mathcal{T}_l \\ \tau \neq k-j \end{matrix}}  t^{\beta_{\tau}-1}E_{\beta_{\tau},\beta_{\tau}}(f_{\tau \tau}(z)t^{\beta_{\tau}})  \Big), & \text{if } j <k.
\end{cases}  
\end{equation}

Similarly, if the matrix symbol $\mathcal{F}(z)$ is upper triangular, then the components of the solution $\mathcal{U}(t)$ take the form
\begin{align}
u_m(t) &= E_{\beta_m}(f_{1 1}(A) t^{\beta_m})\vf_m,\\
u_k(t) &= E_{\beta_k}(f_{k k}(A)t^{\beta_k}) \vf_k  \notag \\
&+ \sum\limits_{j=k+1}^{m} \sum_{l=1}^{m-j-1}  \left[ Q_{k j l}(A) E_{\beta_j}(f_{j j}(A)t^{\beta_j}) \ast \Big(\ast \prod_{\tiny \begin{matrix} \nu \in \mathcal{P}_l \\ \nu \neq j \end{matrix}} t^{\beta_{\nu}-1}E_{\beta_{\nu},\beta_{\nu}}(f_{\nu \nu}(A)t^{\beta_{\nu}})  \Big)\right] \vf_j,
\\
& \hspace{6cm} k=m-1,\dots,1, \notag
\end{align}
where $Q_{k j}(A)$ are operators obtained from $f_{\nu \mu}(A), \nu= k,\dots,m-1, \mu=k+1,\dots,m,$ via the product and sum combinations and $\mathcal{P}_l$ is the set of subsets of the set $\mathcal{P}=\{k, k+1, \dots, m\}$ such that exactly $l$ ($1 \le l \le m-k-1$) numbers except $k$ and $m$ are removed from $\mathcal{P}.$
In turn, the matrix-symbol of the solution operator in representation formula \eqref{case4_S1} 
takes the form
\begin{equation}
s_{kj}(t,z) = 
\begin{cases}
0, & \text{if }  j < k, \\
E_{\beta_k}(f_{k k}(z)t^{\beta_k}), & \text{if } j=k, \\
\sum\limits_{l=1}^{m-j-1} Q_{k j l}(z) E_{\beta_j}(f_{j j}(z)t^{\beta_j}) \ast \Big(\ast \prod\limits_{\tiny \begin{matrix} \nu \in \mathcal{P}_l \\ \nu \neq j \end{matrix}} t^{\beta_{\nu}-1}E_{\beta_{\nu},\beta_{\nu}}(f_{\nu \nu}(z)t^{\beta_{\nu}})  \Big), & \text{if } j > k.
\end{cases}  
\end{equation}


\subsection{Commensurate and non-commensurate rational order systems}
\label{sec_equal}

If all the components of the vector-order $\mathcal{B}$ are equal, then transformation of the matrix-symbol $\mathcal{F}(z)$ to the Jordan canonical form can be effectively utilized in the derivation of representation formulas for the solution. Theorem \ref{thm_gen} serves as an important mathematical base for such an approach. If all the components of $\mathcal{B}$ are rational (not necessarily being equal), then this case can be reduced to the case with equal components, but for the price of increased number of equations (see \cite{Umarov2024}). We note that both these cases were presented in \cite{Umarov2024} for time-fractional systems of pseudo-differential equations. Below, applying the same technique presented in \cite{Umarov2024}, however not going into details, we generalize it for fractional order systems of differential-operator equations, which significantly expands the scope of applications.

Let $\mathcal{F}(A)$ be the matrix-valued operator with the matrix-symbol $\mathcal{F}(z), \ z \in G,$ whose Jordan normal form is
\begin{equation}
\label{J1}
 M^{-1}(z) \mathcal{F}(z) M(z) = \Lambda(z)+N =
 \begin{bmatrix}
 \mathbb{J}_1(z) & 0 & \dots & 0 \\
 0       & \mathbb{J}_2(z) & \dots & 0\\
 \dots & \dots & \dots & \dots \\
 0      &     0   & \dots & \mathbb{J}_L(z)     
 \end{bmatrix},
\end{equation}
where $M(z)$ is a transformation matrix invertible at each $z \in G,$ and $\mathbb{J}_{\ell}(z), \  \ell=1,\dots, L,$ are Jordan blocks corresponding to eigenvalues $\lambda_1(z), \dots, \lambda_L(z)$ of the matrix 
\[
\mathcal{F}(z) = M(z) \Big( \Lambda(z)+N \Big) M^{-1}(z).
\]

First we derive a representation formula for the solution operator of the initial value problem for system \eqref{01} in the homogeneous case
\begin{align} \label{Cauchy_01_h}
D_{\ast}^{\mathcal{B}} {U}(t) &= \mathcal{F}(A) {U} (t), \quad t>0, 
\\
U(0) & =\varPhi, 
\label{Cauchy_02_h}
\end{align}
or due to \eqref{J1} equivalently, 
\begin{align} \label{Cauchy_011_h}
D_{\ast}^{\mathcal{B}} {U}(t) &= M(A) \Big(\Lambda(A) + N\Big) M^{-1}(A) {U} (t), \quad t>0, 
\\
U(0) & =\varPhi, 
\label{Cauchy_021_h}
\end{align}
where $M(A)$  is the transformation matrix-operator with the matrix-symbol $M(z).$ 

In order to solve problem \eqref{Cauchy_011_h},\eqref{Cauchy_021_h} we use the operator approach, considering the following system of ordinary fractional order differential equations dependent on the parameter $z \in G:$
\begin{align} \label{C1_h}
D_{\ast}^{\mathcal{B}} {U}(t,z) &= M(z) \Big(\Lambda(z) + N\Big) M^{-1}(z) {U} (t,z), \quad t>0, \ z \in G,
\\
U(0,z) & =\varPsi, 
\label{C2_h}
\end{align}
assuming that $\varPsi$ is a vector of length $m.$ Since all the components of $\mathcal{B}$ are equal, the matrix-valued operator $ID_{\ast}^{\mathcal{B}} = \mbox{diag} (D_{\ast}^{\beta}, \dots, D_{\ast}^{\beta})$ commutes with $M(z),$ and therefore system \eqref{C1_h} can be expressed as
\begin{equation} 
\label{2C1_h}
M(z) \ ID_{\ast}^{\mathcal{B}} \ M^{-1}(z) {U}(t,z) = M(z) \Big(\Lambda(z) + N\Big) M^{-1}(z) {U} (t,z), \quad t>0, \ z \in G,
\end{equation}
Denote $V(t,z)=M^{-1}(z)U(t,z).$ Then we have the Cauchy problem
\begin{align} 
\label{3C1_h}
D_{\ast}^{\mathcal{B}}  {V}(t,z) &=  \Big(\Lambda(z) + N\Big)  {V} (t,z), \quad t>0, \ z \in G,
\\
V(0,z) & = M(z) \varPsi.
\end{align}
Now applying the Laplace transform in the vector form \eqref{CD_LT_01} to both sides of system \eqref{3C1_h}, we obtain
\begin{align*}
I s^{\mathcal{B}} \mathcal{L}[{V}](s)  &=  I s^{\mathcal{B}-1} M(z) \varPsi +  \Big(\Lambda(z) +N\Big)  \mathcal{L}[V](s), \quad s>0, \ \xi \in G.
\end{align*}
where $I s^{\mathcal{B}}$, $I s^{\mathcal{B}-1}$ are diagonal matrices with diagonal entries $s^{\beta_j}, \ s^{\beta_j-1} \ j=1,\dots,m,$ respectively. 
The latter implies the system of algebraic equations with parameter $z \in G:$
\begin{equation}
\label{n_1}
\Big[ I s^{\mathcal{B}}  - \Big( \Lambda(z) - N \Big) \Big] L[V](s)= I s^{\mathcal{B}-1}  M(z) \varPsi.
\end{equation}
Let 
\begin{equation}
\label{Q0}
Q_0=\{z \in \mathbb{C}: \det \Big[ Is^{\mathcal{B}}  -  \big(\Lambda(z)-N\Big)  \Big]=0, \forall s>0\}.
\end{equation}
If $G \cap Q_0 = \emptyset,$ then system \eqref{n_1} has a unique  solution 
\begin{equation}
\label{s1}
L[V](s) =  \mathcal{N}(s,z) M(z) \varPsi ,
\end{equation}
where 
\begin{align}
{\mathcal{N}}(s,z) &= \Big[ I s^{\mathcal{B}}  -  \Big(\Lambda(z) -N \Big)  \Big]^{-1} I s^{\mathcal{B}-1}.
\label{matrix_03}
\end{align}
This matrix has the block diagonal form 
\[
\mathcal{N}(s,z) =
\begin{bmatrix}
 \mathbb{S}_1(s,z) & 0 & \dots & 0 \\
 0       & \mathbb{S}_2(s,z) & \dots & 0\\
 \dots & \dots & \dots & \dots \\
 0      &     0   & \dots & \mathbb{S}_L(s,z)     
 \end{bmatrix},
 \]
with the blocks $\mathbb{S}_{\ell}(s,z), \ell=1,\dots,L,$ of size $m_{\ell}$ corresponding to the eigenvalue $\lambda_{\ell}$ of multiplicity $m_{\ell}:$
\[
\mathbb{S}_{\ell} =
\begin{bmatrix}
\frac{s^{\beta-1}}{s^{\beta}-\lambda_{\ell}(z)} & \frac{s^{\beta-1}}{(s^{\beta}-\lambda_{\ell}(z))^2} & \dots & \frac{s^{\beta-1}}{(s^{\beta}-\lambda_{\ell}(z))^{m_{\ell}}} \\ 
0  & \frac{s^{\beta-1}}{s^{\beta}-\lambda_{\ell}(z)} & \dots & \frac{s^{\beta-1}}{(s^{\beta}-\lambda_{\ell}(z))^{m_{\ell}-1}} \\
\dots & \dots & \dots & \dots \\
0 & 0 & \dots & \frac{s^{\beta-1}}{s^{\beta}-\lambda_{\ell}(z)}
\end{bmatrix}, \quad \ell=1,\dots L.
\]
Further, applying the inverse Laplace transform, taking into account \eqref{L1} and \eqref{L2}, and returning to $U(t,z)= M(z) V(t,z),$ in accordance with Theorem \ref{thm_gen}, we have 
\begin{equation}
\label{solution_00}
U(t,z) = {E}_{\mathcal{B}}(It^{\mathcal{B}}\mathcal{F}(z)) \varPhi = M(z) \mathbb{E}_{\mathcal{B}} \Big( \big( \Lambda(z) +N \big) t^{\mathcal{B}} \Big) M^{-1}(z)\varPhi, \quad t>0,  
\end{equation}
where $\mathbb{E}_{\mathcal{B}}\Big( \big( \Lambda(z) +N \big) t^{\mathcal{B}} \Big)$ is the block-diagonal matrix of the form
\begin{equation} \label{matrix_04}
\mathbb{E}_{\mathcal{B}}\Big( \big( \Lambda(z) +N \big) t^{\mathcal{B}} \Big) 
= \begin{bmatrix}
	    \mathbb{E}_{\beta}(t^{\beta} J_1(z)) & \dots & 0  \\
	    \dots           & \dots & \dots  \\
	    0  & \dots & \mathbb{E}_{\beta}(t^{\beta} J_L(z))
	    \end{bmatrix},
\end{equation}
with blocks
\begin{align}
&\mathbb{E}_{\beta}(t^{\beta} J_{\ell}(z)) 
=
\begin{bmatrix}
E_{\beta}(\lambda_{\ell}(z)t^{\beta}) & \frac{t^{\beta}E_{\beta}^{\prime}(\lambda_{\ell}(z) t^{\beta})}{1!}  & \dots & \frac{t^{(m_{\ell}-1)\beta}E_{\beta}^{(m_{\ell}-1)}(\lambda_{\ell}(z) t^{\beta})}{(m_{\ell}-1)!}  \\
0 & E_{\beta}(\lambda_{\ell}(z) t^{\beta})  &   \dots   & \frac{t^{(m_{\ell}-2)\beta}E_{\beta}^{(m_{\ell}-2)}(\lambda_{\ell}(z) t^{\beta})}{(m_{\ell}-2)!} \\
\dots & \dots  & \dots  & \dots \\
\dots & \dots & \dots  & \dots \\
0       & \dots &             0                       & E_{\beta}(\lambda_{\ell}(z) t^{\beta})
\end{bmatrix},
\label{m3}
\end{align}
for $\ell=1,\dots,L.$

Thus,  the solution of problem \eqref{Cauchy_01_h}-\eqref{Cauchy_02_h} has the representation
\begin{equation}
\label{solution_0001}
U(t) = S(t,A) \varPhi, \quad t>0, 
\end{equation}
where $S(t,A)$ is the solution matrix operator with the matrix-symbol
\begin{equation}
\label{sol_equal}
\mathcal{S}(t,z) =  {E}_{\mathcal{B}}(It^{\mathcal{B}}\mathcal{F}(z)) = M(z) \ \mathbb{E}_{\mathcal{B}} \Big( t^{\beta} \big( \Lambda(z) +N \big)  \Big) \ M^{-1}(z), \quad t>0, \ z \in G.
\end{equation}

 Now let us
 consider the incommensurate case $\mathcal{B}=(\beta_1,\dots,\beta_m),$ with rational components $\beta_j=q_j/p_j \in (0,1),$ where $p_j, \ q_j$ are positive co-prime integers.  
 Let $p$ be the least common divisor of numbers $p_1,\dots,p_m.$ Then, one can write $\beta_j$ as $\beta_j=n_j/p,$ where $n_j=(q_j p)/p_j$ is integer. Therefore, the operator $D_{\ast}^{\mathcal{B}}$ can be presented in the form
\[
D_{\ast}^{\mathcal{B}}=(D_{\ast}^{\beta_1}, \dots, D_{\ast}^{\beta_m})= \Big((D_{\ast}^{1 \over p})^{n_1}, \dots, (D_{\ast}^{1 \over p})^{n_m} \Big).
\]
It follows that for each $j=1,\dots,m,$ we have
\[
D_{\ast}^{\beta_j} u_j = D_{\ast}^{\frac{1}{p}}  \  \underbrace{\circ \dots \circ}\limits_{n_j \ \text{times}}  \ D_{\ast}^{\frac{1}{p}} u_j.
\]
Introduce a vector-function  $\mathbb{U}(t)$ of length $N=n_1+ \dots + n_m:$ 
\[
\mathbb{U}(t) = (u_1(t), u_1^1(t), \dots, u_1^{n_1-1}(t), \dots, u_m(t), u_m^1(t), \dots, u_m^{n_m-1}),
\]
where $u_j^1 = D_{\ast}^{\frac{1}{p}}u_j, \dots, u_j^{n_j}= D_{\ast}^{\frac{1}{p}}u_j^{n_j-1}, \ j=1,\dots,m. $ 
Now system \eqref{Cauchy_01_h} can be reduced to a system of $N$ equations with equal fractional order $1/p$ derivatives on the left hand side. 
The reduced system has the form
\begin{equation}
\label{augm}
D_{\ast}^{1 \over p} \ \mathbb{U}(t) = \mathbb{F}(A) \mathbb{U}(t),
\end{equation}
where $\mathbb{F}(A)$ is the matrix-operator with the matrix-symbol
\[
\mathbb{F}(z)=
\begin{bmatrix}
\mathbb{F}_{11}(z) & \mathbb{F}_{12}(z)&  \dots & \mathbb{F}_{1m}(z) \\
\mathbb{F}_{21}(z) & \mathbb{F}_{22}(z)&  \dots & \mathbb{F}_{2m}(z) \\
\dots                       &. \dots                     &.  \dots &                        \\
\mathbb{F}_{m1}(z) & \mathbb{F}_{m2}(z)&  \dots & \mathbb{F}_{mm} (z)
\end{bmatrix}
\]
whose diagonal block-matrices are of sizes $n_j \times n_j$:
\[
\mathbb{F}_{jj}(z) = 
\begin{bmatrix}
0 & 1 & 0 & \dots & 0\\
0 & 0 & 1 & \dots & 0\\
\dots & \dots & \dots & \dots & \dots \\
0 & 0 & 0 & \dots & 1\\
f_{jj}(z) & 0 & 0 & \dots & 0
\end{bmatrix}, \quad j=1,\dots,m,
\]
and non-diagonal block-matrices are of sizes $n_j \times n_k$:
\[
\mathbb{F}_{jk}(z) = 
\begin{bmatrix}
0 & 0 & 0 & \dots & 0\\
0 & 0 & 0 & \dots & 0\\
\dots & \dots & \dots & \dots & \dots \\
0 & 0 & 0 & \dots & 0\\
f_{jk}(z) & 0 & 0 & \dots & 0
\end{bmatrix}, \quad j,k=1,\dots, m.
\]
The initial condition for system \eqref{augm} takes the form
\begin{equation}
\mathbb{U}(0) =(\vf_1,  \underbrace{0, \dots, 0}\limits_{n_1-1 \ \text{times}} , \vf_2,\underbrace{0, \dots, 0}\limits_{n_2-1 \ \text{times}}, \dots, \vf_m, \underbrace{0, \dots, 0}\limits_{n_m-1 \ \text{times}}).
\end{equation}

Now, we derive a representation formula for the solution of Cauchy problem \eqref{Cauchy_01_h}, \eqref{Cauchy_02_h}.    
We notice that that the characteristic polynomial of $\mathbb{F}(z),$
\[
P_N(\lambda, z)= \det (I \lambda - \mathbb{F}(z)), \quad \lambda \in \mathbb{C}, \ z \in G,
\] 
and the function $h(s,z)=\det(Is^{\mathcal{B}}-\mathcal{F}(z))$ are connected through the relationship 
\begin{equation}
\label{caseC_1}
h(s,z)=P_N(s^{\frac{1}{p}},z), \quad z \in G.
\end{equation}
Let $\lambda_1(z), \dots, \lambda_L(z)$ be eigenvalues of the matrix $\mathbb{F}(z)$ with respective multiplicities $\mu_1,\dots, \mu_L.$ Then
\[
P_N(\lambda,z)= \prod_{\ell=1}^L (\lambda - \lambda_{\ell}(z))^{\mu_{\ell}},
\]
which implies
\[
\det \Big(I s^{\mathcal{B}}- \mathcal{F}(z)\Big)= \prod_{\ell=1}^L \Big( s^{\frac{1}{p}} - \lambda_{\ell}(z) \Big)^{\mu_{\ell}}.
\]
Similarly, we can write the determinant of the matrix $\mathcal{F}_j(s,z),$ obtained by replacing $j$-th column of the matrix $I s^{\mathcal{B}}-\mathcal{F}(z)$ by the column vector $Is^{B-1} \varPhi,$ in the form
\[
\det(\mathcal{F}_j(s,z))= \sum_{k=1}^m P_{jk}(s^{\frac{1}{p}},z) s^{\beta_k-1} \varphi_k,
\]
where $P_{jk}(\lambda,z)$ is a polynomial in the variable $\lambda$ of order $N-n_j.$ Hence, for the $j$-th component $u_j(t)$ of the Laplace transform of the solution vector $\mathcal{U}(t),$ we have
\begin{align}
L[u_j](s) &= \frac{\det(\mathcal{F}_j(s,z))}{\det \Big(I s^{\mathcal{B}}- \mathcal{F}(z)\Big)} \notag 
\\
&= \sum_{k=1}^m \frac{P_{jk}(s^{\frac{1}{p}},z)}{\prod_{\ell=1}^L \Big( s^{\frac{1}{p}} - \lambda_{\ell}(z) \Big)^{\mu_{\ell}}} s^{\beta_k-1} \varphi_k,
\label{case3_S1}
\end{align}
Further, using the partial fraction decomposition
\[
\frac{P_{jk}(s^{\frac{1}{p}},z)}{\prod_{\ell=1}^L \Big( s^{\frac{1}{p}} - \lambda_{\ell}(z) \Big)^{\mu_{\ell}}} = \sum_{\ell=1}^L \sum_{\nu=1}^{\mu_{\ell}} \frac{C_{\ell \nu}^{jk}(z)}{ \Big( s^{\frac{1}{p}} - \lambda_{\ell}(z) \Big)^{\nu}},
\]
where $C_{\ell \nu}^{j k}(z)$ do not depend on $s,$ equation \eqref{case3_S1} can be expressed as
\[
L[u_j](s) = \sum_{k=1}^m \sum_{\ell=1}^{L} \sum_{\nu=1}^{\mu_L} C_{\ell \nu}^{j k} \frac{s^{\beta_k-1}}{ \Big( s^{\frac{1}{p}} - \lambda_{\ell}(z) \Big)^{\nu}} \ \varphi_k.
\]
Inverting the latter and using formula \eqref{L3}, we obtain
\[
u_j(t)=  \sum_{k=1}^m \sum_{\ell=1}^{L} \sum_{\nu=1}^{\mu_L} \frac{C_{\ell \nu}^{j k} \ t^{\frac{\nu}{p}-\beta_k}}{(\nu-1)!}  E^{(\nu-1)}_{\frac{1}{p} , \frac{1}{p}-\beta_k+1} \Big( t^{\frac{1}{p}}\lambda_{\ell}(z)\Big) \varphi_k, \quad j=1,\dots,m.
\]
It follows that the solution operator  $S(t,A)$ has the matrix-symbol $\mathcal{S}(t,z)$, whose entries are
\begin{equation}
\label{sol_rat}
s_{j k}(t,z) = \sum_{\ell=1}^{L} \sum_{\nu=1}^{\mu_L} \frac{C_{\ell \nu}^{j k} \ t^{\frac{\nu}{p}-\beta_k}}{(\nu-1)!}  E^{(\nu-1)}_{\frac{1}{p} , \frac{1}{p}-\beta_k+1} \Big( t^{\frac{1}{p}}\lambda_{\ell}(z)\Big), \quad j,k=1,\dots,m.
\end{equation}

Summarizing, we obtain the following theorem on formal representations of the solutions:

\begin{thm} 
\label{Repr_2}
Let the eigenvalues $\lambda_{\ell}(\xi), \ \ell=1,\dots, L,$ of the matrix-symbol $\mathcal{F}(\xi)$ of the matrix-valued operator $\mathcal{F}(A)$ in system \eqref{01}, have respective multiplicities $m_{\ell}, \ell =1,\dots,L,$ where $m_1+\dots+ m_{L}=m.$  
Then the solution to system \eqref{01}, \eqref{02} has the representation 
\[
\mathcal{U}(t) = \mathcal{S}(t,A) \Phi + \int_0^{t} \mathcal{S}(t-\tau,A) D_{+}^{1-\mathcal{B}} \mathcal{H}(\tau) d\tau
\]
where $\mathcal{S}(t,A)$ is the matrix-valued solution operator with the matrix-symbol $\mathcal{S}(t,z)$ defined
\begin{enumerate}
\item[(a)]  in \eqref{sol_equal} if $\mathcal{B}=(\beta, \dots, \beta)$ with equal components $\beta \in (0,1],$ and
\item[(b)]  in \eqref{sol_rat} if $\mathcal{B}=(\beta_1,\dots,\beta_m)$ with rational components $\beta_j=q_j/p_j \in (0,1].$\end{enumerate}
\end{thm}

\section{The Riemann-Liouville case}
\label{sec_RL}

Similar results hold in the case when fractional derivatives in system \eqref{01} are in the Riemann-Liouville sense. Therefore, below we briefly formulate the corresponding assertions. 
 
 Consider the initial value problem
 \begin{align}
\label{s215}
&D^{\mathcal B}_+ \mathcal{U}(t) = \mathcal{F} (A) \mathcal{U}(t) + H(t),  \quad t>0,\\
&(\mathcal{J}^{1-\mathcal{B}}\mathcal{U})(0) = \Phi,
\label{s225}
\end{align}  
 where $D^{\mathcal{B}}_+\mathcal{U}(t)=(D_+^{\beta_1}u_1(t),\dots, D_+^{\beta_m}u_m(t)),$ and the matrix-operator $\mathcal{F}(A),$ vector-valued elements $H(t), \ \varPhi$ are specified below.  
 Using the same technique which was used in the case of Caputo derivatives, one can show that in the case of Riemann-Liouville derivatives, the solution matrix-operator $\mathcal{S}_+(t,A)$ has the symbol

 \begin{align}
{\mathcal{S}_+}(t,z) &=
L^{-1}_{s \to t} \left\{  \Big[ I s^{\mathcal{B}}  - \mathcal{F}(z) \Big]^{-1} \right\}, \quad  t \ge 0, \ z \in G \subset \mathbb{C}. 
 \label{matrix_03500} 
\end{align}
 
The following theorems hold: 

\begin{thm}
\label{thm_multi_200}
Let $A$ be a closed operator defined on a Banach space $X,$ the set $G$ satisfies the condition $G \cap Q_0 = \emptyset,$ where $Q_0$ is defined in \eqref{Q0}, $ \varPhi \in \mathrm{Exp}_{A, G},$  and  $H(t) \in AC[\re_+; \mathrm{Exp}_{A, G}].$

Then for any $T>0$ Cauchy problem \eqref{s215}-\eqref{s225}
has a unique solution $U(t)\in C^{\infty}[(0,T]; \mathrm{Exp}_{A, G}] \cap C[[0,T]; \mathrm{Exp}_{A, G}],$ having the representation
\begin{equation}
\label{solution_100}
U(t) = \mathcal{S}_+(t,A)\varPhi + \int\limits_0^t \mathcal{S}_+(t-\tau, A) H(\tau) d\tau, \quad t>0, 
\end{equation}
where $\mathcal{S}_+(t,A)$ is the solution matrix-operator with the matrix-symbol $\mathcal{S}_+(t, z)$ defined in \eqref{matrix_03500}.
\end{thm}

\begin{thm}
\label{thm_multi_200_dual}
Let $X$ be a reflexive Banach space with the conjugate $X^{\ast},$ $A$ be a closed operator with a domain $\mathcal{D} \subset X,$ and $\mathcal{F}(A)$ be a matrix operator with the symbol $\mathcal{F}(z)$ continuous on $G.$ and satisfying condition $G \cup Q_0 = \emptyset,$ where $Q_0$ is defined in \eqref{Q0}. Assume that   $\varPsi \in \mathcal{E}^{\prime}_{A^{\ast}, G^{\ast}},$  and  $H(t) \in AC[\re_+; \mathcal{E}^{\prime}_{A^{\ast}, G^{\ast}}].$

Then for any $T>0$ Cauchy problem
\begin{align} 
\label{Cauchy_01500_dual}
D_{+}^{\mathcal{B}} {V}(t) &= \mathcal{F}(A^{\ast}) {V} (t) + H(t), \quad t>0, 
\\
\mathcal{J}^{1-\mathcal{B}}V(0) & =\varPsi, 
\label{Cauchy_02500_dual}
\end{align}
has a unique solution $V(t) \in C^{\infty}[(0,T]; \mathcal{E}^{\prime}_{A^{\ast}, G^{\ast}}] \cap C[[0,T]; \mathcal{E}^{\prime}_{A^{\ast}, G^{\ast}}],$ having the representation
\begin{equation}
\label{solution_100_dual}
V(t) = \mathcal{S}_+(t,A^{\ast})\varPsi +\int\limits_0^t \mathcal{S}_+(t-\tau, A^{\ast})  H(\tau) d\tau, \quad t>0, 
\end{equation}
where $\mathcal{S}_+(t,A^{\ast})$ is the operator with the matrix-symbol $\mathcal{S}_+(t, z)$ defined in \eqref{matrix_03500}. 
\end{thm}

In what concerns representation formulas for the solution in the case of Riemann-Liouville derivatives, their derivation is similar to the Caputo derivative case. Therefore, we demonstrate the detailed derivation of the representation formula only in the case $m=2$ with arbitrary  multi-order $\mathcal{B}$ and matrix-symbol $\mathcal{F}(z)$, and providing the final result for arbitrary $m \ge 2.$ 

\medskip
\subsection{The case $m=2$ and $\beta_1 \ne \beta_2$}

Consider the system 
\begin{equation}
\label{700}
D_{+}^{\mathcal{B}}\mathcal{U}(t)=\mathbb{F}(A) \mathcal{U}(t) + \mathcal{H}(t),
\end{equation}
where $\mathcal{B}=(\beta_1, \beta_2), 0< \beta_1< \beta_2 \le 1,$ $\mathcal{H}(t)=(h_1(t), h_2(t))$ is a given vector-function, and  
\begin{equation}
\label{701}
\mathbb{F}(A) =
\begin{bmatrix}
f_{11}(A) & f_{1 2}(A) \\
f_{2 1}(A) & f_{2 2}(A)
\end{bmatrix},
\end{equation}
with the initial condition 
\begin{equation}
\label{ic700}
(J^{1-\mathcal{B}} \mathcal{U}) (0)=\varPhi =(\vf_1, \vf_2),
\end{equation} 
where $\varPhi \in  \mathcal{X}_{A,G}.$ We assume that $G$ does not contain roots of the equation 
\[
\Delta(z)={f_{1 1}(z)} f_{2 2}(z) - f_{2 1}(z) f_{1 2}(z)) = 0.
\] 
To find entries of the solution operator $\mathcal{S}(t,z)$ we consider the homogeneous counterpart of system \eqref{200} writing it in the explicit form
\begin{align*}
\begin{cases}
D_{\ast}^{\beta_1} u_1(t) &= f_{1 1}(A) u_1(t) + f_{1 2}(A) u_2(t),\\
D_{\ast}^{\beta_2} u_2(t) & = f_{2 1}(A) u_1(t) + f_{2 2}(A) u_2(t).
\end{cases}
\end{align*}
Applying the Laplace transform and replacing $A$ by the parameter $z$, we have
\begin{equation}
\label{lt700}
\begin{cases}
(Is^{\beta_1}-f_{1 1}(z) ) \mathcal{L}[u_1](s) - f_{1 2}(z) \mathcal{L}[u_2](s) =  \varphi_1 ,\\
- f_{2 1}(z) \mathcal{L}[u_1](s) + (Is^{\beta_2}-f_{2 2}(z) ) \mathcal{L}[u_2](s) =  \varphi_2.
\end{cases}
\end{equation}
The solution of system \eqref{lt200} is
\begin{equation}
\label{rllts21}
\mathcal{L}[u_1](s) = \frac{1}{\Psi(s,z)} \Big(p_1(s,z) \vf_1 + q_{1}(s,z) \vf_2 \Big) \quad z \in G, s > r_{\ast}(z)
\end{equation}
\begin{equation}
\label{rllts22}
\mathcal{L}[u_2](s) = \frac{1}{\Psi(s,z)} \Big( q_{2}(s,z) \vf_1+  p_2(s,z) \vf_2  \Big) , \quad z \in G, s > r_{\ast}(z).
\end{equation}
where $\Psi(s,z)$ is defined in \eqref{lts23} and
\begin{align}
\label{rllts74}
p_1(s,z) &=s^{\beta_2} -f_{2 2}(z), \quad q_1(s,z) = f_{1 2}(z),
\\
\label{rllts75}
p_2(s,z) &=s^{\beta_1} -f_{1 1}(z), \quad q_2(s,z) = f_{2 1}(z),
\end{align}
and $r_{\ast}(z)$ is the real part of the roots of the equation $\Psi(s,z)=0.$ This solution is uniquely defined, since by assumption $G \cap Q = \emptyset,$ where $Q = \{z: \Psi(s,z)=0\}.$

Introduce $\rho_{jk}=j(\beta_2-\beta_1)+(k-j)\beta_2 .$ Obviously, $\rho_{0 0} =0$ and $\rho_{j k} >0$ if $k>0, \ 0 \le j \le k.$
 the entries $s^+_{jl}(t,z), j,l=1,2,$ of  the matrix-symbol $\mathcal{S}_+(t,z)$ have representations:
\begin{align}
s^+_{1 1}(t,z) 
&= \Big(I - f_{2 2}(z)J^{\beta_2}\Big) \sum_{k=0}^{\infty} { \Big( f_{2 2}(z) J^{\beta_2-\beta_1} - \Delta(z)J^{\beta_2} \Big)^k} 
\mathbb{W}_{\beta_1,\beta_1}(t,z),
\label{sol700}
\end{align}

\begin{align}
s^+_{1 2}(t,z) 
&=   f_{1 2}(z)  \Big(I - f_{2 2}(z)J^{\beta_2}\Big) \sum_{k=0}^{\infty}  { \Big( f_{2 2}(z) J^{\beta_2-\beta_1} - \Delta(z)J^{\beta_2} \Big)^k} J^{\beta_2} \mathbb{W}_{\beta_1,\beta_1}(t,z),
\label{sol701}
\end{align}
\begin{align}
s^+_{2 1}(t,z) 
  &= f_{2 1}(z) \Big(I - f_{2 2}(z)J^{\beta_2}\Big) \sum_{k=0}^{\infty}  { \Big( f_{2 2}(z) J^{\beta_2-\beta_1} - \Delta(z)J^{\beta_2} \Big)^k} J^{\beta_2} \mathbb{W}_{\beta_1,\beta_1}(t,z),  
   \label{sol702}
\end{align}
\begin{align}
s^+_{2 2}(t,z) 
&= \sum_{k=0}^{\infty} { \Big( f_{2 2}(z) J^{\beta_2-\beta_1} - \Delta(z)J^{\beta_2} \Big)^k}
 \Big[ \mathbb{W}_{\beta_1,\beta_2}(t,z) 
 -  f_{1 1} (z) J^{\beta_2}
\mathbb{W}_{\beta_1,\beta_1}(t,z)
\Big],
\label{sol703}
\end{align}
where
\[
\mathbb{W}_{\beta_1,\beta_j}(t,z) = \frac{ t^{k \beta_i +\beta_j-1} }{k!} E^{(k)}_{\beta_1, \beta_j} (t^{\beta_1} f_{1 1}(z)), \quad j=1,2.
\]

\begin{thm} \label{thm_m=270}
The solution to system \eqref{s215}, \eqref{s225} has the representation 
\[
\mathcal{U}(t) = \mathcal{S}_+(t,A) \Phi + \int_0^{t} \mathcal{S}_+(t-\tau,A) D_{+}^{1-\mathcal{B}} \mathcal{H}(\tau) d\tau
\]
where $\mathcal{S}_+(t,A)$ is the matrix-valued solution operator with the matrix-symbol $\mathcal{S}_+(t,z),$ entries of which are defined in \eqref{sol700}-\eqref{sol703}.
\end{thm}

\begin{thm}
\label{thm_m=721}
Let $f_{1 2}(z) = 0.$ Then the symbol of the solution operator has entries
\begin{align}
\label{m=711}
s^+_{1 1}(t,z) & =t^{\beta_1-1}E_{\beta_1, \beta_1}(t^{\beta_1}f_{1 1}(z)), \quad S_{1 2}(t,z)=0,  
\\
\label{m=721}
 s^+_{2 1}(t,z) &= f_{2 1}(t,z) \Big( t^{\beta_1-1}E_{\beta_1}(t^{\beta_1}f_{1 1}(z))\Big) \ast \Big( t^{\beta_2-1} E_{\beta_2, \beta_2}(t^{\beta_2}f_{2 2}(z)) \Big),
\\
\label{m=722}
 s^+_{2 2} (t, z) &= t^{\beta_2-1}E_{\beta_2, \beta_2}(t^{\beta_2}f_{2 2}(z)),
 \end{align}
 where $"\ast"$ is the convolution operation.
\end{thm}
{\it Proof} is similar to the proof of Theorem \ref{thm_m=21}.  

\medskip
\subsection{The case $m \ge 2$ and arbitrary $\mathcal{B}$}

In the case of arbitrary $m \ge 2,$ following the same approach demonstrated in Section \ref{sec_3.2}, one can derive the representation formula. Namely, instead of equation \eqref{am03} one has
\begin{equation}
\label{rlam03}
(Is^{\mathcal{B}} - \mathcal{F}(z)) \mathcal{L} [ \mathcal{U}] (s) =  \varPhi,
\end{equation}
and therefore, instead of equation \eqref{SOLUTION_SYMBOL} one obtains
\begin{equation}
\label{RLSOLUTION_SYMBOL}
\mathcal{S}(t,z) = \mathcal{L}^{-1}_{s \to t} \Big[ (Is^{\mathcal{B}} - \mathcal{F}(z))^{-1} \Big],
\end{equation}
the symbol of the solution operator. The entries of the latter are
\begin{align}
s_{j l}(t,z) &= \sum_{k=0}^{\infty} \sum\limits_{\alpha} R_{\alpha, j, l, k}(z) \mathcal{L}_{s \to t}^{-1} \left[ \frac{1}{s^{\nu_{\alpha, j, l, k}}\Big(s^{\beta_1}+g(z)\Big)^{k+1}}\right] \notag
\\
&= \sum_{k=0}^{\infty} \sum\limits_{\alpha} R_{\alpha, j, l, k}(z) J^{\nu_{\alpha, j, l, k}} \left[ \frac{t^{k\beta_1 +\beta_1-1}}{k!}  E^{(k)}_{\beta_1, \beta_1} \Big(-t^{\beta_1}g(z) \Big) \right], \quad j, l =1,\dots,m,
\label{rlsolm5}
\end{align}
where $R_{\alpha, j, l, k}(z)$ are the sum and product combinations of entries of $\mathcal{F}(z).$

\subsection{The case $m \ge 2$ and $\mathcal{B}=(\beta,\dots,\beta)$}

If the components of $\mathcal{B}$ are equal, then again we can use Jordan normal form to derive a representation formula for the solution. In this case, following the method used in Section \ref{sec_equal}, we consider the system of ordinary fractional order differential equations depending on the parameter $z \in G$
\begin{equation}
\label{B5}
D_{+}^{\mathcal{B}}  {V}(t,z) = \Big( \Lambda(z) +N \Big) {V} (t,z), \quad t>0, 
\end{equation}
where the symbol $\Lambda(z) + N$ of the matrix operator $\Lambda(A)+ N$ has representation in the Jordan block form \eqref{J1}.  
Then for the Laplace transform of $V(t)$ we obtain a linear system of algebraic equations

\begin{equation}
\label{n_25}
\Big(I s^{\mathcal{B}} - \Lambda(z) - N\Big) L[V](s)=   M^{-1}\varPhi.
\end{equation}
If $Q_0 \cap G = \emptyset,$ then the latter has a unique solution represented through the inverse matrix $\mathcal{N}(s,z) = \Big(I s^{\mathcal{B}} - \Lambda(z) - N\Big)^{-1},$ which has the block diagonal form 

\[
\mathbb{S}(s,z) =
\begin{bmatrix}
 \mathbb{S}_1(s,z) & 0 & \dots & 0 \\
 0       & \mathbb{S}_2(s,z) & \dots & 0\\
 \dots & \dots & \dots & \dots \\
 0      &     0   & \dots & \mathbb{S}_L(s,z)     
 \end{bmatrix},
 \]
with the blocks $\mathbb{S}_{\ell}(s,z), \ell=1,\dots,L,$ of size $m_{\ell}$ corresponding to the eigenvalue $\lambda_{\ell}$ of multiplicity $m_{\ell}:$
\[
\mathbb{S}_{\ell} =
\begin{bmatrix}
\frac{1}{s^{\beta}-\lambda_{\ell}(z)} & \frac{1}{(s^{\beta}-\lambda_{\ell}(z))^2} & \dots & \frac{1}{(s^{\beta}-\lambda_{\ell}(z))^{m_{\ell}}} \\ 
\, & \, & \, & \\
0  & \frac{1}{s^{\beta}-\lambda_{\ell}(z)} & \dots & \frac{1}{(s^{\beta}-\lambda_{\ell}(z))^{m_{\ell}-1}} \\
\dots & \dots & \dots & \dots \\
0 & 0 & \dots & \frac{1}{s^{\beta}-\lambda_{\ell}(z)}
\end{bmatrix}, \quad \ell=1,\dots L.
\]
Now using formula \eqref{L4} one can find the inverse Laplace transform of each entries of matrices $\mathbb{S}_{\ell}, \ \ell=1,\dots,L.$ Hence, the solution matrix-operator has the block-matrix representation 
\[
\mathcal{S}(t,A)= M(A) \mathbb{G}_{\mathcal{B}}\Big( \big( \Lambda(A) +N \big) t^{\mathcal{B}} \Big) M^{-1}(A),
\]
where the block matrix-operator $\mathbb{G}_{ \mathcal{B}}\Big( \big( \Lambda(A) +N \big) t^{\mathcal{B}} \Big) $ has the matrix symbol
\begin{equation} \label{matrix_045}
\mathbb{G}_{ \mathcal{B}}\Big( \big( \Lambda(z) +N \big) t^{\mathcal{B}} \Big) 
= \begin{bmatrix}
	    \mathbb{G}_{\beta}(t^{\beta} J_1(z)) & \dots & 0  \\
	    \dots           & \dots & \dots  \\
	    0  & \dots & \mathbb{G}_{\beta}(t^{\beta} J_L(z))
	    \end{bmatrix},
\end{equation}
with blocks
\begin{align}
&\mathbb{G}_{\beta}(t^{\beta} J_{\ell}(z)) 
=
\begin{bmatrix}
t^{\beta-1}E_{\beta, \beta}(\lambda_{\ell}(z)t^{\beta}) & \frac{t^{2\beta-1}E_{\beta, \beta}^{\prime}(\lambda_{\ell}(z) t^{\beta})}{1!}  & \dots & \frac{t^{(m_{\ell}-1)\beta-1}E_{\beta, \beta}^{(m_{\ell}-1)}(\lambda_{\ell}(z) t^{\beta})}{(m_{\ell}-1)!}  \\
0 & t^{\beta-1}E_{\beta, \beta}(\lambda_{\ell}(z) t^{\beta})  &   \dots   & \frac{t^{(m_{\ell}-2)\beta-1}E_{\beta, \beta}^{(m_{\ell}-2)}(\lambda_{\ell}(z) t^{\beta})}{(m_{\ell}-2)!} \\
\dots & \dots  & \dots  & \dots \\
\dots & \dots & \dots  & \dots \\
0       & \dots &             0                       & t^{\beta-1}E_{\beta, \beta}(\lambda_{\ell}(z) t^{\beta})
\end{bmatrix},
\label{rlm3}
\end{align}
for $\ell=1,\dots,L.$

Concluding,  the solution matrix-operator 
has the representation
\begin{equation}
\label{rlsolution_0001}
U(t) = S(t,A) \varPhi, \quad t>0, 
\end{equation}
where $S(t,A)$ is the solution matrix operator with the matrix-symbol
\begin{equation}
\label{sol_0001}
\mathcal{S}(t,z) =  M(z) \ \mathbb{G}_{\mathcal{B}} \Big( \big( \Lambda(z) +N\big) t^{\mathcal{B}} \Big) \ M^{-1}(z), \quad t>0, \ z \in G.
\end{equation}

\section{Some applications and examples}
\label{sec_appl}

\begin{ex} {Time-fractional systems of ordinary differential equations.} \end{ex}
Consider the following initial-value problem for time-fractional system of ordinary differential equations
\begin{align*}
D_{\ast}^{\mathcal{B}} U(t) &= \mathbb{A} \, U(t) + H(t), \quad t>0,
\\
U(0) &= U_0,
\end{align*}
where $\mathbb{A}$ is a (constant) $m \times m$-matrix, $U_0=(u^0_1, \dots, u^0_m) \in \re^m,$ and $H(t)$ is an absolute continuous vector-function. Theorems obtained above are applicable to this case with the proper interpretation. 

Let $\mathcal{B}=(\beta_1, \dots, \beta_m), \ 0< \beta_j \le 1,$ be arbitrary numbers,  $a_{j k}, \ j, k =1, \dots, m,$  be entries of the matrix $\mathbb{A},$  and $d=\det(\mathbb{A}).$ Suppose $\beta_1 = \min \{\beta_1,\dots,\beta_m\}.$ Consider first the corresponding homogeneous system
\[
D_{\ast}^{\mathcal{B}} U(t) = \mathbb{A} \, U(t), \quad U(0) = U_0.
\]
Define the function
\[
\Psi(s) = \det (I s^{\beta} - \mathbb{F}) = s^{\beta_1+\dots + \beta_m} + g s^{\beta_2 +\dots + \beta_m } + \dots.
\]

The solution to the latter has the form $U(t)=S(t) U_0,$ where $S(t), t \ge 0,$ is the solution matrix. The components of $S(t)$ due to Theorem \ref{thm_001} has entries 
\begin{align}
s_{j l}(t) 
&= \sum_{k=0}^{\infty} \sum\limits_{\alpha} Q_{\alpha, j, l, k} J^{\nu_{\alpha, j, l, k}} \left[ \frac{t^{k\beta_1}}{k!}  E^{(k)}_{\beta_1} \Big(-g t^{\beta_1} \Big) \right], \quad j, l =1,\dots,m,
\label{solm500}
\end{align}
where $Q_{\alpha, j, l, k}$ are defined similar to $Q_{\alpha, j, l, k}(z)$ in \eqref{q}, replacing $f_{kj}(z)$ by $a_{k j}.$ In particular, if $m=2,$ then
\begin{align*}
s_{1 1}(t)
&= \Big(I - a_{2 2}J^{\beta_2}\Big)\sum_{k=0}^{\infty} { \Big( a_{2 2} J^{\beta_2-\beta_1} - d J^{\beta_2} \Big)^k} \Big[\frac{ t^{k \beta_1}}{k!}  E^{(k)}_{\beta_1}(a_{1 1} t^{\beta_1} ) \Big],
\end{align*}
\begin{align*}
s_{1 2}(t) &=   a_{1 2} \sum_{k=0}^{\infty}  { \Big( a_{2 2} J^{\beta_2-\beta_1} - dJ^{\beta_2} \Big)^k}  J \Big( \frac{t^{k\beta_1+\beta_1-1}}{k!}   E^{(k)}_{\beta_1, \beta_1} (  a_{1 1}t^{\beta_1} ) \Big),
\end{align*}
\begin{align*}
s_{2 1}(t) &= a_{2 1}\sum_{k=0}^{\infty} { \Big( a_{2 2} J^{\beta_2-\beta_1} - dJ^{\beta_2} \Big)^k}
  J^{\beta_2} \Big(\frac{ t^{k \beta_1}}{k!}   E^{(k)}_{\beta_1} ( a_{1 1} t^{\beta_1} ) \Big),
\end{align*}
\begin{align*}
s_{2 2}(t) &= \sum_{k=0}^{\infty} { \Big( a_{2 2} J^{\beta_2-\beta_1} - d J^{\beta_2} \Big)^k}
 \Big[ \frac{ t^{k \beta_1}}{k!}   E^{(k)}_{\beta_1} (  a_{1 1}t^{\beta_1} )
-  a_{1 1}  J \Big(\frac{ t^{k \beta_1+\beta_1-1}}{k!}    E^{(k)}_{\beta_1, \beta_1} (t^{\beta_1} a_{1 1}) \Big) \Big] .
\end{align*}
Additionally, if $a_{12}=0,$ then
\[
s_{1 1}(t)=E_{\beta_1} (a_{1 1}t^{\beta_1}), \quad s_{1 2}(t) = 0, 
\]
\[
s_{2 1}(t) = a_{2 1} E_{\beta_1} (a_{1 1}t^{\beta_1}) \ast (t^{\beta_2-1} E_{\beta_1, \beta_2} (a_22 t^{\beta_2})), \quad s_{2 2} =E_{\beta_2} (a_{2 2}t^{\beta_2}).
\]

\begin{ex} Blood alcohol  level problem.\end{ex} 
The authors of paper \cite{Alcohol} considered the following blood alcohol problem using fractional order derivatives in the sense of Caputo-Djrbashian:
\begin{equation}
\begin{cases}
D^{\alpha}_{\ast} A(t) = -\kappa_1 A(t),\\
D^{\beta}_{\ast} B(t) = \kappa_1 A(t) -\kappa_2 B(t),
\end{cases}
\end{equation}
where $A$ represents the concentration of alcohol in the stomach and $B$ is the
concentration of alcohol in the blood, and $\kappa_1, \ \kappa_2$ are some real constants, which indicate transition rates. The initial conditions are given by
\[
A(0)=A_0, \ B(0)=B_0. \footnote{$B_0=0$ if initially there is no alcohol in the blood.}
\]
This problem can be presented as $D_{\ast}^{\mathcal{B}}\mathcal{U}(t)=\mathcal{F} \mathcal{U}(t),$ where $\mathcal{U}(t)=(A(t), B(t)),$ $\mathcal{B}=(\alpha, \beta)$ and 
\[
\mathcal{F} = 
\begin{bmatrix}
-\kappa_1 & 0 \\
\kappa_1 & - \kappa_2
\end{bmatrix}.
\]

In accordance with Theorem \ref{thm_m=20} the solution $\mathcal{U}(t)$ has the representation
\begin{align}
\label{alcoh1}
 A(t) &= A_0 E_{\alpha}(-\kappa_1 t^{\alpha}),
\\
\label{alcoh2}
 B(t) &= \kappa_1 A_0  \Big( E_{\alpha}(- \kappa_1 t^{\alpha}) \Big) \ast \Big( t^{\beta-1}E_{\beta,\beta} (-\kappa_2 t^{\beta}) \Big) + B_0 E_{\beta}(-\kappa_2 t^{\beta}).
\end{align}

We note that authors of \cite{Alcohol} found the solution in the form $A(t)=A_0 E_{\alpha}(-\kappa_1 t^{\alpha})$ and (with $B_0=0$)
\[
B(t) = \beta \int_0^t \kappa_1 A_0 E_{\alpha}(-\kappa_1 (t-s)^{\alpha}) s^{\beta-1}E_{\beta}^{\prime}(-\kappa_2 s^{\beta}) ds,
\]
which is the same as \eqref{alcoh1}, \eqref{alcoh2}, due to the equality $\beta E_{\beta}^{\prime}( z) = E_{\beta,\beta}(z).$

\medskip

\begin{ex} Fractional order systems for a relativistically free particle. \end{ex}
The wave function of a relativistically free particle of mass $m$ is described by the Klein-Gordon equation
\[
\left( \frac{1}{c^2} \frac{\partial^2}{\partial t^2} - \nabla^2 + \frac{m^2 c^2}{\hbar^2} \right) \Psi (t,x) = 0, \ x \in \re^3,
\]
where $c$ is the speed of light, $\hbar$ is Planck's constant, and $\nabla = (\frac{\partial}{\partial x_1}, \frac{\partial}{\partial x_2}, \frac{\partial}{\partial x_3}),$ the gradient operator. Dirac's equation for relativistically free particle, in fact, is a system of the form 
\[
i\hbar \frac{\partial \Psi (x,t)}{\partial t} = (i c \hbar {\bf \alpha} \cdot \nabla + \beta m c^2) \Psi(x,t),
\]
where ${\bf \alpha}$ and $\beta$ are $4 \times 4$ matrices satisfying certain conditions, and $\Psi(x,t)$ is a multi-component wave function. Using the adopt units $c=\hbar=1,$ the latter can be written in the explicit form \cite{CollasKlein}
\begin{equation}
\label{Dirac_3}
i\begin{bmatrix}
\frac{\partial \psi_1}{\partial t} \\
\frac{\partial \psi_2}{\partial t} \\
\frac{\partial \psi_3}{\partial t} \\ 
\frac{\partial \psi_4}{\partial t}
\end{bmatrix}
=
\begin{bmatrix}
m  & 0 & -i \frac{\partial }{\partial x_3} & i\frac{\partial}{\partial x_1} + \frac{\partial }{\partial x_2}  \\
0 & m & - i \frac{\partial }{\partial x_1} + \frac{\partial}{\partial x_2} & \frac{\partial }{\partial x_3} \\
 -i \frac{\partial}{\partial x_3} & - i\frac{\partial}{\partial x_1} - \frac{\partial }{\partial x_2} & - m & 0 \\
 -i \frac{\partial}{\partial x_1} + \frac{\partial}{\partial x_2}  & i \frac{\partial}{\partial x_3} & 0 &-m  
\end{bmatrix}
\begin{bmatrix}
\psi_1 \\
\psi_2 \\
\psi_3 \\ 
\psi_4
\end{bmatrix}
\end{equation}
Multiplying by $-i,$ we can rewrite system \eqref{Dirac_3} in the form \eqref{01}:
\begin{equation}
\label{Dirac_1}
\frac{\partial \Psi (t,x)}{\partial t} = \mathcal{F}(A) \Psi(t,x), \quad t>0, \ x \in \re^3,
\end{equation}
where $\mathcal{F}(A), \ A= -i \nabla = -i (\frac{\partial}{\partial x_1}, \frac{\partial}{\partial x_2}, \frac{\partial}{\partial x_3}),$ is the matrix-valued operator with the symbol
\begin{equation}
\label{Dirac_2}
\mathcal{F}(\xi)=
 \begin{bmatrix}
-i m  & 0 &  \xi_3 & i\xi_1 + \xi_2  \\
0 & -i m & - i \xi_1+\xi_2 & \xi_3 \\
 -i \xi_3 & - i\xi_1 -\xi_2 & i m & 0 \\
 -i \xi_1 + \xi_2  & i \xi_3 & 0 &i m  
\end{bmatrix}, \quad \xi=(\xi_1,\xi_2,\xi_3) \in \re^3.
\end{equation}

Replacing $\partial/\partial t$ on the left  of \eqref{Dirac_1} by $D^{\alpha}_{\ast}$ we obtain Dirac-like fractional order system. Note that some Dirac-like systems are considered in papers \cite{TP,Ercan,VazquezMendes2003,FSKG}.

Thus, consider the system 
\[
D^{\alpha}_{\ast} \Psi(t,x)= \mathcal{F}(-i \nabla) \Psi(t,x),
\] 
where $1/2 < \alpha \le 1$ and  
$\mathcal{F}(-i \nabla)$ has the symbol $\mathcal{F}(\xi_2,\xi_2,\xi_3)$ defined in \eqref{Dirac_2}. 
The matrix $\mathcal{F}(z)$ has eigenvalues 
\[
\lambda_{1-4}(\xi) = \pm \sqrt{-m^2 \pm i \sqrt{(\xi_1^2+\xi_2^2)^2 + \xi_3^4}}.
\]
If $\xi \neq 0,$  then all the eigenvalues are of multiplicity one, and hence diagonalizable. Consequently, there exists an invertible matrix $M(\xi), $ such that  $\mathcal{F}(\xi)=M(\xi) \Lambda(\xi) M^{-1}(\xi),$ with $\Lambda (\xi)=\mbox{diag}(\lambda_1(\xi), \dots, \lambda_4(\xi)).$
Thus, in accordance with Theorem \ref{Repr_2} the solution is represented in the form  $\Psi(t,x)=\mathcal{S}(t, -i \nabla) \Psi(0, x),$ where the solution pseudo-differential operator $\mathcal{S}(t,-i\nabla)$ has the matrix-symbol
\[
\mathcal{S}(t,\xi)= M(\xi) \mbox{diag} \Big(E_{\alpha} (t^\alpha) \lambda_1(\xi)), \dots, E_{\alpha} (t^{\alpha}\lambda_4(\xi))  \Big) M^{-1}(\xi), 
\]
and components $\psi_k(x), k=1,\dots, 4,$ of $\Psi(0,x)$ have Fourier transforms with compact supports in $G = \re^3 \setminus \{0\}.$

\medskip

\begin{ex} A commensurate system of pseudo-differential equations.  \end{ex}

Let the matrix-valued operator $\mathcal{F}(A)$  
is given by the matrix-symbol 
\[
\mathcal{F}(z) = 
\begin{bmatrix} 
-|z|^2 & a^2 \bar{z} \\
z      & -|z|^2
\end{bmatrix}, \quad z \neq 0.
\]
Consider the system with matrix-valued operator on the right corresponding to the symbol $\mathcal{F}(z):$ 
\begin{equation}
\label{ex11}
D_{\ast}^{\mathcal{B}}\mathcal{U}(t)=
\begin{bmatrix} 
-A^2 & a^2 A^{\ast} \\
A     & -A^2
\end{bmatrix} 
\mathcal{U}(t), \quad t>0, 
\end{equation}
with $\mathcal{B}=(\beta,\beta), \ 0< \beta \le 1,$ $A^{\ast}$ is the adjoint of $A,$ and the initial condition 
\begin{equation}
\label{ex2}
\mathcal{U}(0)=(\vf_1,\vf_2).
\end{equation}
Assume that $a \neq 0.$ Then, one can easily verify that eigenvalues of $\mathcal{F}(z)$ are $\lambda(z) = |z|^2 \pm a|z|.$ Therefore, 
$
\mathcal{F}(z)= M(z) \Lambda (z)  M^{-1} (z),
$
where
\[
\Lambda(z) = 
\begin{bmatrix} 
|z|^2 - a |z| & 0 \\
0     & |z|^2 + a |z|
\end{bmatrix},
\quad 
M(z) = 
\begin{bmatrix} 
a\bar{z} & a\bar{z} \\
-|z|     & |z|
\end{bmatrix}.
\]
Then due to Theorem \ref{Repr_2} the solution operator $\mathcal{S}(t,A)$ has the matrix symbol
\begin{align*}
\mathcal{S}(t,z) &= 
\begin{bmatrix} 
a\bar{z} & a\bar{z} \\
-|z|     & |z|
\end{bmatrix}
\begin{bmatrix} 
E_{\beta}((-|z|^2-a|z|) t^{\beta-1}) & 0 \\
0    & E_{\beta}((-|z|^2 +a|z|) t^{\beta-1})
\end{bmatrix}
\begin{bmatrix} 
 \frac{1}{2a\bar{z}} &  -\frac{1}{2|z|} \\
 \frac{1}{2a \bar{z}}     & \frac{1}{2|z|}
\end{bmatrix} \\
\, \\
&= 
\begin{bmatrix} 
\frac{1}{2} \Big(E^{-}(z,t) + E^+(z,t)\Big) &  -\frac{a \bar{z}}{2|z|} \Big( E^{-}(z,t) - E^+(z,t) \Big) \\
-\frac{|z|}{2a\bar{z}} \Big( E^{-}(z,t) - E^+(z,t) \Big) & \frac{1}{2} \Big(E^{-}(z,t) +\frac{1}{2} E^+(z,t) \Big)
\end{bmatrix},
\end{align*}
where $E^{\pm}(z,t)=E_{\beta} \Big((-|z|^2 \pm a |z|) t^{\beta-1} \Big).$

Suppose $A=\sqrt{-\Delta},$ where $\Delta$ is the Laplace operator with the domain 
\[
Dom(A)= H^{1}(\re^n) = \{f \in L_2(\re^n): \int\limits_{\re^n} (1+|\xi|^2)^{1/2} |F[f](\xi)|^2 d\xi < \infty \}.
\] 
Here $F[f]$ is the Fourier transform of $f(x).$ It is known that the spectrum of $A$ is the positive semi-axis, and hence we can accept $G=[0, \infty)$. Then, the symbol $\mathcal{S}(t,z)$ simplifies to
\begin{align*}
\mathcal{S}(t,z) &= 
\begin{bmatrix} 
\frac{1}{2} \Big(E^{-}(z,t) + E^+(z,t)\Big) &  -\frac{a}{2} \Big( E^{-}(z,t) - E^+(z,t) \Big) \\
-\frac{1}{2a} \Big( E^{-}(z,t) - E^+(z,t) \Big) & \frac{1}{2} \Big(E^{-}(z,t) +\frac{1}{2} E^+(z,t) \Big)
\end{bmatrix},
\end{align*}

If $a = 0,$ then $\mathcal{F}(z)$ has a double eigenvalue $\lambda(z) = |z|^2$ and has the Jordan form
$
\mathcal{F}(z)= M(z) \Big(\Lambda (z) + N \Big) M^{-1} (z).
$
where
\[
\Lambda(z) = 
\begin{bmatrix} 
|z|^2 & 0 \\
0     & |z|^2
\end{bmatrix},
\quad 
N(z) = 
\begin{bmatrix} 
0 & 1 \\
0     & 0
\end{bmatrix},
\quad
M(z) = 
\begin{bmatrix} 
0 & 1/z \\
1     & 1
\end{bmatrix}.
\]
Assume for simplicity that $\mathcal{B}=(1/2, 1/2).$  
Then, in accordance with Theorem \ref{Repr_2}, the solution operator $S(t,A)$ has the matrix-symbol 
\begin{align*}
\mathcal{S}(t,z) &= 
\begin{bmatrix} 
0 & 1/z \\
1     & 1
\end{bmatrix}
\begin{bmatrix} 
E_{1/2}(-|z|^2 t^{1/2}) & t^{1/2} E_{1/2}^{\prime}(-|z|^2 t^{1/2}) \\
0    & E_{1/2}(-|z|^2 t^{1/2})
\end{bmatrix}
\begin{bmatrix} 
 -{z} &  1 \\
 {z}     & 0
\end{bmatrix} \\
&= 
\begin{bmatrix} 
E_{1/2}(-|z|^2 t^{1/2}) &  0 \\
z t^{1/2}E^{\prime}_{1/2}(-|z|^2 t^{1/2}) & E_{1/2}(-|z|^2 t^{1/2})
\end{bmatrix}. 
\end{align*}
Thus, the solution $\mathcal{U}(t)$ has components
\begin{align}
\label{exsol1}
u_1(t) &= E_{1/2}(-A^2 t^{1/2}) \vf_1,\\
u_2(t) &= t^{1/2}  A E_{1/2}^{\prime}(A^2 t^{1/2}) \vf_1+E_{1/2}(-A^2 t^{1/2}) \vf_2.
\label{exsol2}
\end{align}

\medskip

\begin{ex} An incommensurate system of pseudo-differential equations. \end{ex}
Now assume $\mathcal{B}=( 1/2, 1/3 ),$ $A=\sqrt{-\Delta},$ 
and $\mathcal{F}(z)$ is as in Example 4 with $a=0.$ Then, since $\beta_1=1/2$ and $\beta_2=1/3$ are rational numbers, we can use the technique described in Section \ref{sec_equal} and reduce problem \eqref{ex1}, \eqref{ex2} to a system of five  equations with equal orders $\beta_j^{\ast} = 1/6, j=1,\dots,5.$  The reduced system has the form 
\begin{equation}
\label{ex3}
\begin{bmatrix} 
D_{\ast}^{1/6} U_1(t,x)  \\
D_{\ast}^{1/6} U_2(t,x) \\
D_{\ast}^{1/6} U_3(t,x) \\
D_{\ast}^{1/6} U_4(t,x) \\
D_{\ast}^{1/6} U_5(t,x) \\
\end{bmatrix}
=
\begin{bmatrix} 
0 & 1 & 0 & 0 & 0 \\
0 & 0 & 1 & 0 & 0 \\
\Delta & 0 & 0 & 0 & 0 \\
0 & 0 & 0 & 0 & 1 \\
\sqrt{-\Delta} & 1 & 0 & \Delta & 0 \\
\end{bmatrix} 
\begin{bmatrix} 
 U_1(t,x)  \\
U_2(t,x) \\
 U_3(t,x) \\
 U_4(t,x) \\
U_5(t,x) \\
\end{bmatrix}, \quad t>0, \ x \in \re^n,
\end{equation}
with the initial condition
\[
\Big( U_1(0,x),  U_2(0,x),  U_3(0,x) , U_4(0,x),  U_5(0,x) \Big) = \Big( \vf_1(x),0,0,\vf_2(x),0\Big).
\]
The components $U_1(t,x) = u_1(t,x)$ and $U_4(t,x)=u_2(t,x)$ correspond to the solution of \eqref{ex1}, \eqref{ex2}. Applying the Fourier and Laplace transforms we can transform system \eqref{ex3} to the following algebraic equations
\begin{align*}
s^{1/6} L[F[U_1]](s,z) &= L[F[U_2]](s,z) + s^{-5/6} F[\vf_1](z), \\
s^{1/6} L[F[U_2]](s,z) &= L[F[U_3]](s,z), \\
s^{1/6} L[F[U_3]](s,z) &= - |z|^2 L[F[U_1]](s,z), \\
s^{1/6} L[F[U_4]](s,z) &= L[F[U_5]](s,z) + s^{-5/6} F[\vf_2](z), \\
s^{1/6} L[F[U_5]](s,z) &=  z L[F[U_1]](s,z) - |z|^2 L[ F[U_4]](s,z).
\end{align*}
It follows that 
\[
L[F[U_1]] (s,z) = \frac{s^{-1/2} F[\vf_1](z)}{s^{1/2}+|z|^2},
\]
and
\[ 
L[F[U_4]] (s,z) = z \frac{s^{-1/2} F[\vf_1](z)}{(s^{1/2}+|z|^2)(s^{1/3}+|z|^2)} + \frac{s^{-2/3} F[\vf_2](z)}{s^{1/3}+|z|^2}.
\]
Thus the solution $\mathcal{U}(t,x)$ of problem \eqref{ex1}, \eqref{ex2} in accordance with Theorem \ref{Repr_2} has the following components:
\begin{align*}
u_1(t,x) &= E_{1/2}(\Delta t^{1/2}) \vf_1(x),\\
u_2(t,x) &=  \sqrt{-\Delta} \Big( E_{1/2}(\Delta t^{1/2})\ast t^{-2/3}E_{1/3,1/3} (\Delta t^{1/3}) \Big) \vf_1(x)+E_{1/3}(\Delta t^{1/3}) \vf_2(x).
\end{align*}

Note that Theorem \ref{thm_m=21} is applicable for this problem, as well, resulting in the same solution.



\begin{thebibliography}{999} \rm


\bibitem{DalKrein} S.G. Krein.
Linear differential equations in Banach space. AMS, Providence,
R.I., 1971.


\bibitem{VasPis}
V.V. Vasil'yev, S.I. Piskarev, Differential equations in Banach
spaces I. Theory of cosine operator functions II. J. Math. Sci. 122,
2 (2004), 3055-3174.

\bibitem{SKM93}
{S.G. Samko, A.A. Kilbas, O.I. Marichev}, {Fractional Integrals and
Derivatives: Theory and Applications},  Gordon and Breach Science
Publishers, New York - London, 1993.

\bibitem{Ba2001}
E. Bazhlekova, Fractional evolution equations in Banach spaces,
Dissertation, Technische Universiteit Eindhoven, 2001.

\bibitem{Kos93}
{V.A. Kostin},  The Cauchy problem for an abstract differential
equation with fractional derivatives, {Russ. Dokl. Math.} {46}
(1993) 316-319.

\bibitem{Koc89}
{A.N. Kochubei}, A Cauchy problem for evolution equations of
fractional order, {Differential Equations} {25} (1989) 967-974.

\bibitem{GLZ99}
{R. Gorenflo, Yu. F. Luchko, P.P. Zabreiko}, On solvability of
linear fractional differential equations in Banach spaces. { Fract.
Calc. Appl. Anal.}, {2} (1999) 163-176.

\bibitem{EidKoch2004}
S. D. Eidelman, A. N. Kochubei, Cauchy problem for fractional
diffusion equations, Journal of Differential Equations, 199 (2004)
211--255.

\bibitem{KilbasST}
A.A.Kilbas, H.M. Srivastawa, J.J.Trijillo, Theory and applications
of fractional differential equations, Elsevier Science, 2006.

\bibitem{Umarov_book2015} Umarov, S.: Introduction to Fractional and Pseudo-Differential Equations with Singular Symbols. Springer (2015)

\bibitem{Veber}  Veber, V.K.: The structure of general solution of the system $y(\alpha) = Ay, 0 < \alpha \le 1.$
Trudy Kirgiz. Gos. Univ. Ser. Mat. Nauk, 1976, 11, 26–32 (in Russian).

\bibitem{Garrapa} Garrappa, R., Popolizio, M.: Computing the matrix Mittag-Leffler function with applications to fractional Calculus. J. Sci. Comput. {\bf 77}, 129--153 (2018). https://doi.org/10.1007/s10915-018-0699-5

\bibitem{Varsha} Daftardar-Gejji, V., Babakhani, A.: Analysis of a system of fractional differential equations. 
{J. of Math. Anal. and Appl} {\bf 293}, 511--522 (2004)

\bibitem{UACh} Umarov, S., Ashurov, R. Chen, Y. On a Method of Solution of Systems of Fractional Pseudo-Differential Equations. Fract Calc Appl Anal 24, 254–277 (2021). https://doi.org/10.1515/fca-2021-0011

\bibitem{DengLiGuo} Deng, W.,  Li, Ch.,  Guo, Q.: Analysis of fractional differential equations with multi-orders.  
{Fractals} {\bf 15}(2), 173--182 (2007)

\bibitem{good} Odibat, Z.: Analytic study on linear systems of fractional differential equations. Computers and Mathematics with Applications {\bf 59}, 1170-1183 (2010)

\bibitem{Umarov2024} Umarov, S.: Representations of solutions of systems of time-fractional pseudo-differential equations, Frac. Calc. Appl. Anal. (2024) (to appear)

\bibitem{DasGupta} Das, S., Gupta, P.K.: A mathematical model on fractional Lotka-Volterra equations. 
{Journal of Theoretical Biology} {\bf 277}(1), 1--6 (2011)


\bibitem{Rihan} Rihan, F.:
Numerical modeling of fractional-order biological systems. 
{Abstract and Applied Analysis} {\bf 2013}, 1--13 (2013)


\bibitem{GuoFang} Guo, Ch., Fang, Sh.: Stability and approximate analytic solutions of the fractional Lotka-Volterra equations for three competitors. 
{Advanced Difference Equations} {\bf 219}, 1--14 (2016)

\bibitem{Khan} Khan, N.A., Razzaq, O. A., Mondal, S.P., Rubbab, Q.:
Fractional order ecological system for complexities of interacting species with harvesting threshold in imprecise environment. 
{Advances in Difference Equations} {\bf 405}, 1--34 (2019)


\bibitem{Rana} Rana, S., Bhattacharya, S., Pal, J., Guerekata, G., Chattopadhyay, J.: Paradox of enrichment: A fractional differential approach with memory. 
{Physica A: Statistical Mechanics and its Applications} {\bf 392}(17), 3610--3621 (2013)


\bibitem{Zeb} Zeb, A., Zaman, G., Chohan, M.I., Momani, Sh., Erturk, V.S.:
Analytic numeric solution for SIRC epidemic model in fractional order.  
{Asian J. of Math. and Appl.} {\bf 2013}, 1--19 (2013)


\bibitem{Islam} Islam, R., Pease, A., Medina, D., Oraby, T.: 
Integer versus fractional order SEIR deterministic and stochastic models of measles. 
{International Journal of Environmental Research and Public Health} {\bf 17}(6), 1--19 (2020)

\bibitem{TP} Pierantozzi, T.: Fractional evolution Dirac-like equations: Some properties and a discrete Von Neumann-type analysis. Journal of Computational and Applied Analysis {\bf 224}(1), 284--295 (2009)


\bibitem{Ercan} Ercan, A.: On the fractional Dirac systems with non-singular operators. Thermal Science {\bf 23} (6), 2159--2168, 2019

\bibitem{Bernstein} Bernstein, S.: Fractional Dirac operator. Noncommutative Analysis, Operator Theory and Applications {\bf 252}  27–41, 2016.



\bibitem{GKMR}  Gorenflo, R.,  Kilbas, A.A.,  Mainardi, F.,  Rogosin S.V.: Mittag-Leffler Functions, Related Topics and Applications. Springer 2nd ed. (2020)

\bibitem{ML} Haubold, H.J.,  Mathai, A.M., Saxena, R.K.: Mittag-Leffler functions and their applications.
{Journal of Applied Mathematics} {\bf 2011}, 1--51 (2011). doi: 10.1155/2011/298628

\bibitem{DSh} Dunford, N., Schwartz, J.: Linear Operators, Interscience Publishers, 1964.


\bibitem{Robertson}
A.P. Robertson, W.J. Robertson, Topological Vector Spaces, Cambridge
University Press, London, 1964.

\bibitem{Radino}
Ya.V.Radyno, Linear equations and bornology, BSU, Minsk, 1982 (in
Russian)

\bibitem{Um98}
{S.R. Umarov},  Nonlocal boundary value problems for pseudo-differential and differential operator equations II, {Differential Equations} {34} (1998) 374-381.

\bibitem{Sato} Kaneko A.: On the structure of hyperfunctions with compact supports. Proc. Japan Acad. {\bf 47}, 956--958, 1971.


\bibitem{Umarov_handbook} Umarov, S: Fractional Duhamel principle. In "Handbook of Fractional Calculus", De-Gruyter, 2019.


\bibitem{Alcohol} Almeida, R., Bastos, N.R.O., Monteiro, M.T.T.:
Modeling some real phenomena by fractional differential equations
Math Methods Appl Sci, 39 (16) (2016), pp. 4846-4855, doi: 10.1002/mma.3818

\bibitem{CollasKlein} Collas, P.,  Klein, D.: 
The Dirac Equation in Curved Spacetime: A Guide for Calculations, Springer, 2019.

\bibitem{VazquezMendes2003} 
Vazquez, L.,  Mendes, R.V.: Fractionally coupled solutions of the diffusion equation. 
{Appl. Math. Comp.} {\bf 141}, 125--130 (2003)

\bibitem{FSKG} Blackledge, J., Babajanov B.: The fractional Schr\"odinger-Klein-Gordon equation and intermediate relativism.
Mathematica Aeterna, {\bf 3} (8),  601--615, 2013.

\end{thebibliography}
\end{document}